\documentclass[pdflatex,sn-mathphys-num]{sn-jnl}

\usepackage{graphicx}%
\usepackage{multirow}%
\usepackage{amsmath,amssymb,amsfonts}%
\usepackage{amsthm}%
\usepackage{mathrsfs}%
\usepackage[title]{appendix}%
\usepackage{xcolor}%
\usepackage{textcomp}%
\usepackage{manyfoot}%
\usepackage{booktabs}%
\usepackage{algorithm}%
\usepackage{algorithmicx}%
\usepackage{algpseudocode}%
\usepackage{listings}%
\usepackage{pgf,tikz}
\usetikzlibrary{arrows}

\theoremstyle{thmstyleone}%
\newtheorem{theorem}{Theorem}
\newtheorem{proposition}[theorem]{Proposition}%

\theoremstyle{thmstyletwo}%
\newtheorem{example}{Example}%
\newtheorem{remark}{Remark}%
\newtheorem{lemma}[theorem]{Lemma}
\theoremstyle{thmstylethree}%
\newtheorem{definition}{Definition}%
\newtheorem{corollary}[theorem]{Corollary}
\begin{document}

\title[Simple and subdirectly irreducible weakly dicomplemented lattices]{ Simple,  subdirectly irreducible weakly dicomplemented lattices}


\author*[1]{\fnm{Yannick Lea} \sur{Jeufack Tenkeu}}\email{yannick.tenkeu@facsciences-uy1.cm}

\author*[2]{\fnm{Leonard} \sur{Kwuida}}\email{leonard.kwuida@bfh.ch}
\equalcont{These authors contributed equally to this work.}


\affil[1]{\orgdiv{University of Yaound\'{e} 1}, \orgname{Department of Mathematics}, \orgaddress{P.O. BOX 812}, 
\city{Yaound\'{e}}, 
\country{Cameroon}}

\affil[2]{\orgdiv{
	Bern University of Applied Sciences}, \orgname{School of Business}, \orgaddress{\street{Brückenstrasse 73}, \city{Bern}, \postcode{3005}, \country{Switzerland}}}

\abstract{In this work, we exhibit 
several subclasses of weakly dicomplemented lattices (WDLs) based on their skeletons and dual skeletons. We investigate normal filters (resp. ideals) and show that the set of normal filters (resp. ideals) forms a complete lattice, which is not a sublattice of the lattice of all filters (ideals).  The normal filter (ideal) generated by a subset and the join of two normal filters (resp. ieals) are characterized. We further prove that the lattice of normal filters is isomorphic to the lattice of normal ideals, and that the only class of filters (or ideals) that generate a congruence in WDLs is the class of normal filters. For distributive WDLs, the congruences generated by filters are characterized. Using normal filters, we characterize simple, subdirectly irreducible, and regular WDLs.
	Moreover, it is shown that the congruences generated by normal filters are permutable, and that regular distributive WDLs are congruence-permutable and verify the congruence extension property (CEP). Finally, we prove that, under certain conditions, the lattice of normal filters is isomorphic to the lattice of filters of the Boolean center of a distributive WDL. It is also established that the lattice of normal filters of a WDL $L$
	embeds into the lattice of normal filters of the power $L^{X}$ of $L$.
		}

\maketitle

\section{Introduction}\label{sec1}

Dicomplemented lattices originate from Formal Concept Analysis \cite{Rudolf5} and were introduced by Rudolf Wille 
to capture the equational theory of concept algebras. Weakly dicomplemented lattices (WDLs, for short) were introduced in \cite{kwdd} from the equations that generate all equations that are up-to-date known to hold in all concept algebras. WDLs generalize Boolean algebras and distributive double pseudocomplemented lattices. While concept algebras are complete lattices, WDLs are not necessarily complete  lattices\cite{kwdd}.
From a categorical point of view, it is important to study substructures, products, and morphisms in order to gain a deeper understanding of the structure of concept algebras arising from a given context, with possible applications in computer science and data analysis.\par 

In ordered algebraic structures, congruences, filters, ideals, and homomorphisms are fundamental tools that have been extensively used to study these algebras and have been particularly intensively investigated by many researchers. Gr$\ddot{a}$tzer, in \cite{gratz2, gratz}, investigated lattice congruences and characterized those generated by filters in distributive lattices, proving that each filter in a distributive lattice corresponds to the class of the unit element of some congruence. This characterization provides further results concerning the decomposition of distributive lattices. R. Beazer and M.E. Adams in \cite{MADAMRBEASER, RBEAZER}, characterized congruences generated by normal filters in distributive double 
$p$-algebras and deduced characterizations of simple and subdirectly irreducible distributive double $p$-algebras.
Tenkeu et al. in \cite{GaelT2,Tenk}	characterized simple and subdirectly irreducible double Boolean algebras using congruences generated by particular filters and ideals in double Boolean algebras. Kalmbach, in \cite{kalmbach}, investigated and established several interesting properties of orthomodular lattices using congruences generated by filters.

From the preceding results, it is evident that congruences are the main tools for the decomposition of algebraic structures. Therefore, a good description of congruences constitutes a substantial contribution to understanding the structural theory of WDLs and can aid in directly characterizing simple and subdirectly irreducible concept algebras.

In the case of WDLs, several algebraic investigations have been carried out. In \cite{kwdd}, Kwuida proved the prime ideal theorem for WDLs and introduced several classes of filters (and ideals), such as primary and normal filters, thereby initiating the study of congruences in WDLs, particularly for distributive concept algebras.

In a bounded distributive lattice, every filter 
$F$ is the cokernel of at least one congruence; that is, there exists a lattice congruence 
$\theta$ such that
$F=\{x\in L\mid (x,1)\in\theta \}$, denoted by 
$[1]_{\theta}$. However, this property does not generally hold in WDLs. Thus, it becomes important to characterize those filters of WDLs that are congruence cokernels, to describe the congruences generated by these filters, and to use this characterization for a direct description of regular, simple, and subdirectly irreducible distributive WDLs.
\par 
In the present work, we focus on the investigation of certain subclasses of WDLs, following the properties of their skeletons and dual skeletons, normal filters, congruences generated by filters, as well as regular, simple, and subdirectly irreducible WDLs. In this framework, we also examine essential aspects of normal filters, such as normal filters generated by subsets, the join of two normal filters, and the relationship between normal filters and normal ideals. This study leads to the following contributions:
\begin{itemize}
	\item[(1)] The construction of an extensive list of examples of WDLs and the introduction of several subclasses of WDLs.
	\item[(2)] A study of certain properties of normal filters and the relationships between normal filters and normal ideals.
	\item[(3)] An analysis of congruences generated by filters and some related properties.
	\item[(4)] The characterization of simple, regular, and subdirectly irreducible distributive WDLs.
\end{itemize}
This paper is organized into five sections. Section 2 presents the basic background notions on WDLs and some known results from lattice theory. In Section 3, new subclasses of WDLs are introduced, and a list of examples is provided. Section 4 establishes that the lattice of normal filters is not a sublattice of the lattice of filters and is isomorphic to the lattice of normal ideals of a given WDL. In this section, the normal filter generated by a subset and the join of two normal filters are also characterized. Section 5 proves that normal filters are the only class of filters corresponding to the class of unity of some congruence in a WDL. Congruences generated by normal filters are then characterized, and some related properties are derived, including the fact that two induced congruences commute. Finally, the lattice of normal filters is used to characterize simple and subdirectly irreducible distributive WDLs. The paper concludes with an investigation of some properties of the quotient of a WDL generated by normal filters and show that the class of distributive regular WDLs has the congruence extension property.

\section{Preliminaries}\label{sec2}

	The main results of this paper concern the characterization of congruences generated by filters, as well as the notions of simplicity and subdirect irreducibility in the variety of weakly dicomplemented lattices. We begin by establishing some general notations. It is assumed that the reader is familiar with equivalence relations, order relations, lattices, distributive lattices, and Boolean algebras. For further information on these concepts from universal algebra, we refer the reader to \cite{Burris}.

We denote by $\underline{A}$ an algebra $(A; F)$, where $A$ (the \textbf{universe}) is a nonempty set and $F$ (the \textbf{set of fundamental operations}) is a set of finitary operations on $A$. We let $\Delta_{A} := \{(x, x) \mid x \in A\}$ and $\nabla_{A} := A \times A$ denote, respectively, the smallest and the largest equivalence relations on $A$.
A congruence relation on $\underline{A}$ is an equivalence relation on $A$ that is compatible with the fundamental operations of $\underline{A}$. We denote by $\operatorname{Con}(\underline{A})$ the set of all congruences of $\underline{A}$. Ordered by inclusion, this set forms a complete lattice with smallest element $\Delta_{A}$ and largest element $\nabla_{A}$.	We denote by $\underset{i \in I}{\Pi}A_{i}$ the direct product of the algebras $(\underline{A}_{i})_{i \in I}$. The algebra $\underline{A}$ is said to be (directly) indecomposable if it is not isomorphic to a direct product of algebras, each having at least two elements.

Birkhoff introduced the notion of a subdirect product as a powerful generalization of the direct product. An algebra $\underline{A}$ is a subdirect product of $(\underline{A}_{i})_{i \in I}$ if:
(a) $\underline{A}$ is a subalgebra of $\underset{i \in I}{\Pi}\underline{A}_{i}$, and
(b) $\pi_{i}(\underline{A}) = \underline{A}_{i}$ for each $i \in I$,
where $\pi_{i}$ denotes the $i$-th projection. A map $\alpha : \underline{A} \to \underset{i \in I}{\Pi}\underline{A}_{i}$ is called a subdirect embedding if $\alpha(\underline{A})$ is a subdirect product of $(\underline{A}_{i})_{i \in I}$. The algebra $\underline{A}$ is subdirectly irreducible if, for every subdirect embedding $\alpha : \underline{A} \to \underset{i \in I}{\Pi}\underline{A}_{i}$, there exists some $i \in I$ such that the composition $\pi_{i} \circ \alpha: \underline{A} \to \underline{A}_{i}$ is an isomorphism (see \cite{Burris}, p. 57).
The importance of such algebras is highlighted by a classical theorem of Birkhoff (see \cite{Burris}, p. 58), which states that in any equational class of algebras, every algebra can be embedded into a direct product of subdirectly irreducible algebras. A particular case of a subdirectly irreducible algebra is a simple algebra, namely one whose lattice of congruences forms a two-elements chain $\{\Delta_{A}, \nabla_{A}\}$.
There is a strong connection between the decomposition of an algebra and the existence of certain congruences on it. In the rest of the paper, we assume that 
$A$ has at least two elements. In this case, a convenient characterization of subdirect irreducibility was given by Garrett Birkhoff.  An algebra 
$\underline{A}$	is subdirectly irreducible when  $Con(\underline{A})\setminus \{\Delta_{A}\}$  has a least element (\cite{Burris}, Theorem 8.4, p.57), called the monolith of 
$Con(\underline{A})$.
 We conclude this section by recalling some basic definitions and results on weakly dicomplemented lattices
 (WDLs) and Formal Concept Analysis, presenting concept algebras as a large subclass of WDLs.
\begin{definition}
	A \textbf{formal context} $\mathbb{K}:=(G, M, I)$ consists of two sets $G$ and $M$ and a relation $I\subseteq G\times M$. The elements of $G$ are called the \textbf{objects}, and the elements of $M$ are called \textbf{attributes} of the context.  To express that  an object $g$ is in relation $I$ with an attribut $m$, we write $gIm$ or $(g, m)\in I$ and read it as'' the object $g$ has the attribute $m$".\\
	For a set $A\subseteq G$ of objects,  we define $A':=\{m\in M\mid gIm~\text{for all }g\in A\},$ the set of attributes common to the objects in $A$. Correspondingly, for a set $B\subseteq M$ of attributes, we define 
	$B':=\{g\in G\mid gIm~\text{for all}~ m\in B\},$
	the set of objects having attributes in $B$.\par 
	A \textbf{formal concept} of a context $\mathbb{K}:=(G, M, I)$ is a pair $(A, B)$ with $A\subseteq G$ and $B\subseteq M$ such that $A'=B$ and $B'=A$.  For a formal concept $(A,B)$, we define:
	\begin{itemize}\item   
		its \textbf{weak negation} by $(A,B)^{\Delta}:=((G\setminus A)'', (G\setminus A)'),$
		\item  its \textbf{weak opposition} by 
		$ (A,B)^{\nabla}:=((M\setminus B)',(M\setminus B)''),$
		\item nothing $0:=(M', M),$
		\item all  $1:=(G, G')$.
	\end{itemize}
	The algebra	$\mathcal{A}(\mathbb{K}):=(\mathfrak{B}(\mathbb{K}); \wedge,\vee,^{\Delta},^{\nabla}, 0,1)$ is called the \textbf{concept algebra} of the context $\mathbb{K}$, where $\mathfrak{B}(\mathbb{K})$ is the set of formal concepts of $\mathbb{K}$,  $\wedge$ and $\vee$ denote the  infimum and  supremum operations of the concept lattice defined by $$(A, B)\vee (C, D)=((A\cup C)^{''}, B\cap D)\text{ and}~ (A, B)\wedge (C, D)=(A\cap C, (B\cup D)'').$$
\end{definition}
\begin{definition}\cite{kwdd}\label{theo:concept} 
	A \textbf{weakly dicomplemented lattice} is an algebra $(L; \wedge,\vee,^{\Delta},^{\nabla},0,1)$ of type $(2,2,1,1,0,0)$ for which $(L;\wedge,\vee,0,1)$ is a bounded lattice satisfying the following  properties:
	\vspace{0.2cm}\\
	\begin{minipage}{7cm}\begin{itemize}
			\item[(1)]  $x^{\Delta\Delta}\leq x$,
			\item[(2)] $x\leq y\implies y^{\Delta}\leq x^{\Delta}$,
			\item[(3)] $(x\wedge y)\vee(x\wedge y^{\Delta})=x$,
		\end{itemize}
	\end{minipage}\begin{minipage}{6cm}
		\begin{itemize}	\item[(1')] $x^{\nabla\nabla}\geq x$,
			\item[(2')] $x\leq y\implies y^{\nabla}\leq x^{\nabla}$,
			\item[(3')] $(x\vee y)\wedge (x\vee y^{\nabla})=x$.
		\end{itemize}
	\end{minipage}~\\
	
	We call $x^{\Delta}$ (resp. $x^{\nabla}$)  the \textbf{weak } (resp. \textbf{dual weak) complement} of $x$. The operations $^{\Delta}, ^{\nabla}$ and 
	$(^{\Delta}, ^{\nabla})$ are called \textbf{weak complementation}, \textbf{dual weak complementation} and   \textbf{weak dicomplementation}.
	We recall some important concepts and ordered algebraic structures related to weakly dicomplemented lattices. 
\end{definition}	

\begin{itemize}
	\item Let $(L, \leq)$ be a partially ordered set. A \textbf{closure operator} $c$ on $L$ is a map ${c: (L, \leq)} \to (L, \leq)$ such that for all $x, y \in L$, the following three conditions hold: \\
	(i) $x \leq c(x)$, \quad (ii) $x \leq y$ implies $c(x) \leq c(y)$, \quad and \quad (iii) $c(c(x)) = c(x)$.
	
	\item If we replace condition (i) by (i)' $c(x) \leq x$, then $c$ is termed an \textbf{interior operator} on $L$.
	
	\item (\cite{Sylvia-Karl})\label{def:ortho} An \textbf{orthocomplemented lattice} (or \textbf{ortholattice} for short) is a bounded lattice with a unary operation $^{\bot}: L \to L$ such that, for all $a, b \in L$:
	$$~(i)~a \leq b \implies b^{\bot} \leq a^{\bot},~
	(ii)~(a^{\bot})^{\bot} = a,~ 
	(iii)~a \vee a^{\bot} = 1~ \text{and}~ 
	(iv) ~a \wedge a^{\bot} = 0.$$	
	\item An ortholattice becomes an \textbf{orthomodular lattice} (OML) if the orthomodular law holds:  
	$$\text{(OML)} \quad a \leq b \Rightarrow b = a \vee (a^{\bot} \wedge b).$$
	
	\item Let $L$ be a bounded lattice and $x\in L$. The element $x^{\ast}\in L$ (resp. $x^{+}\in L$) is the peseudocomplement (resp. dual pseudocomplement) of $x$ if 
	$$x\wedge y=0 \Leftrightarrow y\leq x^{\ast}  ( \text{resp.}~ x\vee y=1~\Leftrightarrow y\geq x^{+}), \text{for all}~x\in L.$$
	\item (\cite{RBEAZER}) A \textbf{pseudocomplemented lattice} (resp. \textbf{dual pseudocomplemented lattice}), also called a $p$-algebra (resp. \textbf{dual p-algebra}), is an algebra $(L; \wedge, \vee, ^{\ast}, 0, 1)$ (resp. $(L; \wedge, \vee, ^{+}, 0,1)$) of type $(2,2,1,0,0)$ in which every element $x$ has a pseudocomplement ( resp. dual pseudocomplement) that is $x^{\ast}$ (resp. $x^{+}$). 
	When the lattice is distributive, the $p$-algebra (resp. dual p-algebra)  is said to be \textbf{distributive}.
	
	\item A \textbf{double $p$-algebra} is an algebra $(L; \wedge, \vee, ^{\ast}, ^{+}, 0, 1)$ of type $(2,2,1,1,0,0)$ such that $(L; \wedge, \vee, ^{\ast}, 0, 1)$ is a $p$-algebra and $(L; \wedge, \vee, ^{+}, 0, 1)$ is a \textbf{dual $p$-algebra} (\cite{T.Katrinak5,T.Katrinak2}
\end{itemize}
It is proved that every distributive double $p$-algebra is a  WDL, but not all double $p$-algebras are WDLs (see Example 1.1.3 in \cite{kwdd}).  

The map $x \mapsto x^{\nabla\nabla}$ (resp. $x \mapsto x^{\Delta\Delta}$) defines a closure operator (resp. interior operator) on $\textbf{L}$.  

The set 
$$S(\textbf{L}) := \{x \in L \mid x^{\nabla\nabla} = x\} \quad \text{(resp. } \overline{S}(\textbf{L}) := \{x \in L \mid x^{\Delta\Delta} = x\}\text{)}$$ 
are called the \textbf{skeleton} (resp. \textbf{dual skeleton}) of $\textbf{L}$, and the set 
$$B(L) := \{x \in L \mid x^{\Delta} = x^{\nabla}\}$$ 
is called the subset of elements with negation, or the \textbf{center} of $L$.  

One defines the operations $\sqcap$, $\sqcup$, $\overline{\sqcap}$, and $\underline{\sqcup}$ on $L$ by:
\[
x \sqcap y := (x^{\nabla} \vee y^{\nabla})^{\nabla}, \quad
x \sqcup y := (x^{\nabla} \wedge y^{\nabla})^{\nabla}, \quad
x \overline{\sqcap} y := (x^{\Delta} \vee y^{\Delta})^{\Delta}, \quad \text{and}~
x \underline{\sqcup} y := (x^{\Delta} \wedge y^{\Delta})^{\Delta}.
\]

Let   $\mathcal{L}_{1} = (L_{1}, \leq_{1})$ and $\mathcal{L}_{2} = (L_{2}, \leq_{2})$ be two lattices.  The lattices  $\mathcal{L}_{1}$ and $\mathcal{L}_{2}$ are order  \textbf{isomorphic} if there exists a map $\phi: L_{1} \to L_{2}$ that is  onto and satisfies the equivalence  
$$a \leq_{1} b  \iff \phi(a) \leq_{2} \phi(b), a, b\in L_{1}.$$
Let $(L; \wedge, \vee, ^{\Delta}, ^{\nabla}, 0,1)$ be a WDL, we denote by   $Con(L)$ its set of congruences. For $\theta \in Con(L)$ and $a \in L$, we write
$[a]_{\theta} := \{x \in L \mid (a, x) \in \theta\}$
for the equivalence class containing $a$.

\begin{proposition}(\cite{kwdd}, Prop. 1.3.4)\label{prop:ortho} 
	For any WDL $\textbf{L}$, the following statements hold:
	\begin{itemize}
		\item[(1)] The algebras 
		$(S(\textbf{L}); \wedge, \sqcup, ^{\nabla}, 0, 1) \quad \text{and} \quad 
		(\overline{S}(\textbf{L}); \overline{\sqcap}, \vee, ^{\Delta}, 0, 1)$ 
		are ortholattices.
		\item[(2)] The algebra 
		$(B(L); \vee, \wedge, ', 0, 1)~\text{with}~ x':=x^{\Delta}=: x^{\nabla}$ 
		is a Boolean algebra,  called \textbf{Boolean center} of $L$ (\cite{kwdd}, Cor. 3.3.8).
	\end{itemize}
\end{proposition}
From now on, $L$ stands for a WDL $(L; \wedge, \vee, ^{\Delta}, ^{\nabla}, 0,1)$. Some basic properties of WDLs are collected in the following proposition.  

\begin{proposition} \cite{kwdd}\label{prop: axiom} 
	Let $x, y \in L$. The following statements hold: \vspace{0.3cm}
	
	\begin{minipage}{6cm}
		\begin{itemize}
			\item[(4)] $(x \wedge y)^{\Delta} = x^{\Delta} \vee y^{\Delta}$.
			\item[(5)] $x \overline{\sqcap} y = (x \wedge y)^{\Delta\Delta}$. 
			\item[(6)] $x^{\Delta} \leq y \Leftrightarrow y^{\Delta} \leq x$.
			\item[(7)] $y \wedge x = 0 \Rightarrow y^{\Delta} \geq x$.
			\item[(8)] $x \vee x^{\Delta} = 1$.
			\item[(9)] $x^{\nabla} \leq x^{\Delta}$.
		\end{itemize}
	\end{minipage}
	\begin{minipage}{7cm}
		\begin{itemize}
			\item[(4')] $(x \vee y)^{\nabla} = x^{\nabla} \wedge y^{\nabla}$.
			\item[(5')] $x \sqcup y = (x \vee y)^{\nabla\nabla}$.
			\item[(6')] $x^{\nabla} \geq y \Leftrightarrow y^{\nabla} \geq x$.
			\item[(7')] $y \vee x = 1 \Rightarrow y^{\nabla} \leq x$.
			\item[(8')] $x^{\nabla} \wedge x = 0$.
		\end{itemize}
	\end{minipage}
\end{proposition}
It is easy to show that	$L$ satisfies the identity 
$(\star):~ x = x^{\Delta\nabla}$ 
if and only if $L=B(L)$.
A nonempty subset $F$ of $L$ is called an \textbf{order filter}  if $x\in F, x\leq y\Rightarrow y\in F, \forall x, y\in L$. An \textbf{order ideal} is defined dually.  $F$ is a \textbf{filter} if $F$ is an order filter such that  for all $x, y \in L$:  
$x, y \in F \Rightarrow x \wedge y \in F$. 	An \textbf{ideal} is defined dually.  	
A \textbf{proper filter} (resp. \textbf{ideal}) is a filter (resp. ideal) that is different from $L$. It is well known that the set of filters (resp. ideals) of a WDL, ordered by set inclusion, forms a complete lattice, denoted by $F(L)$ (resp. $I(L)$). The join of two filters (resp. ideals) is calculated as in the case of lattices (see \cite{gratz2}), Boolean algebra (\cite{Burris}, Lemma 3.9).

\section{Some subclasses of WDLs}~

In order to investigate simple and subdirectly irreducible WDLs, we first provide several examples of WDLs, examine some specific subclasses of weakly dicomplemented lattices, and explore some of their properties.

\subsection{Some examples of  WDLs}~\\
 The following known result gives a usefull method to construct weakly dicomplemented lattice in a given bounded lattice. For an element $v$ of a complete lattice $\mathcal{L}:=(L, \leq)$ we define 
$$v_{\ast}:=\bigvee\{x\in L\mid x< v\}~\text{and}~v^{\ast}:=\bigwedge\{x\in L\mid v< x \}.$$ We call $v$ $\bigvee$-irreducible (or supremum irreducible) (resp. $\bigwedge$- irreducible or infimum irreducible) if $v\neq v_{\ast}$ (resp. $v\neq v^{\ast}$). $J(L)$ denotes the set of all $\bigvee$-irreducible elements and $M(L)$ the set of all $\bigwedge$-irreducible elements. A subset $X\subseteq L$ is called \textbf{supremum- dense} in $L$ if every element from $L$ can be  represented as a supremum of the subset $X$, and dually, \textbf{infimum-dense} if $v=\bigwedge\{x\in X\mid v\leq x\}$ (\cite{GanterVille}, pp.7).
\begin{proposition}
	(\cite{Kwuidahajime} Proposition 3) \label{propp:wdlv} Let $L$ be a finite lattice.  Define two unary operations $^{\Delta}$ and $^{\nabla}$ on $L$ by 
	$$x^{\Delta}:=\bigvee\{a\in J(L)\mid a\nleq x\}~\text{and}~x^{\nabla}:=\bigwedge\{m\in M(L)\mid x\nleq m\}.$$
	Then $(L; \wedge, \vee, ^{\Delta}, ^{\nabla}, 0, 1)$ is a WDL. In general, for $J(L)\subseteq G$ and $M(L)\subseteq H$, the operations $^{\Delta_{G}}$ and $x^{\nabla_{H}}$ defined by 
	$$x^{\Delta_{G}}:=\bigvee\{a\in G\mid a\nleq x\}~\text{and}~^{\nabla_{H}}:=\bigwedge\{m\in H\mid x\nleq m\}$$ turn $(L; \wedge, \vee, ^{\Delta_{G}}, x^{\nabla_{H}}, 0, 1)$ into a WDL. 
\end{proposition} 

\begin{example}
	\begin{enumerate}
		
		\item  (\cite{kwdd}) For a Boolean algbera $(B; \wedge, \vee, ', 0,1)$, the algebra $(B; \wedge, \vee, ',', 0,1)$ (complementation duplicated, i.e $x^{\Delta}:=x'=: x^{\nabla})$ is a WDL.\item The algebra $\mathcal{A}(\mathbb{K})$ of concepts of a given context $\mathbb{K}$ is a WDL (\cite{Rudolf5}).
		
		\item There exists a unique WDL structure on an $n$-elements chain, $n \geq 2$. The trivial weak complement (resp. dual weak complement) of $x \in L$ is defined by $x^{\Delta} = 1$ if $x \neq 1$, and $x^{\Delta} = 0$ otherwise (resp. $x^{\nabla} = 0$ if $x \neq 0$, and $0^{\nabla} = 1$).
		
			\item Consider the bounded lattice $\mathcal{M}_{4,2}$ in Fig. 1. One can check that $J(M_{4,2})=\{a, b,1\}$ and $M(M_{4,2})=\{e,a,b\}$. Using Proposition \ref{pro: irre}, we   define  		
		the unary operations $^{\Delta}$ and $^{\nabla}$  on $M_{4,2} = \{0, a, b, e, 1\}$ by:  
		$x^{\Delta} = 1, \ x \in \{0, a, b, e\}, \quad 1^{\Delta} = 0,$ 
		and  
		$0^{\nabla} = 1, \quad x^{\nabla} = 0 \ \text{for} \ x = e, 1, \quad a^{\nabla} = b, \quad b^{\nabla} = a$. By  Proposition \ref{propp:wdlv}, the algebra $\mathcal{M}_{4,2} = (M_{4,2}; \vee, \wedge, ^{\Delta}, ^{\nabla}, 0,1)$ is a distributive WDL.
		Moreover, 
		$S(M_{4,2}) = \{0, a, b, 1\}$ and $\overline{S}(M_{4,2}) = \{0, 1\}$, both of which are Boolean algebras.
		
		\item Consider the two elements Boolean algebra $\mathcal{C}_{2}$ and the three elements chain $\mathcal{C}_{3}$ endowed with the unique structure of WDL. Then  $\mathcal{P}_{6}=\mathcal{C}_{2}\times \mathcal{C}_{3}$  with $P_{6} = \{0, u, v, a, b, 1\}$, represented in Fig. 2. is a WDL as product of WDLs. 
		The unary operations $^{\Delta}$ and $^{\nabla}$ are defined by the Cayley table given in (T:2). One has   
		$S(\mathcal{P}_{6}) = \overline{S}(\mathcal{P}_{6}) = \{0, u, b, 1\}$, which forms  a Boolean algebra.

		\item Consider the lattice $L_{6}$ represented in Fig. 2, with  
		$L_{6} = \{0, u, v, a, b, 1\}$, and the unary operations $^{\Delta}$ and $^{\nabla}$ defined by the Cayley table given in (T:1).  
		The algebra $(L_{6}; \wedge, \vee, ^{\Delta}, ^{\nabla}, 0,1)$ is a WDL, with  
		$S(L_{6}) = \{0, u, v, 1\}$ and $\overline{S}(L_{6}) = \{0, u, b, 1\}$.  
		Both $S(L_{6})$ and $\overline{S}(L_{6})$ are Boolean algebras, and their intersection is the three-elements chain $C_{3} = \{0, u, 1\}$. Moreover, $L_{6} \neq S(L_{6}) \cup \overline{S}(L_{6})$.
		\item Consider the lattice $\mathcal{K}_{7}$,  represented in Fig.3. with $K_{7} = \{0, u, v, a, b, w, 1\}$.  One can show that $J(K_{7})=\{w,u,v\}$ and $M(K_{7})=\{u, v, a, b\}$. Using Proposition  \ref{propp:wdlv} one defines the operation $^{\Delta}$ and $^{\nabla}$ on $K_{7}$ by $u^{\nabla}=v, v^{\nabla}=u, x^{\nabla}=0, x\in\{a, b, w, 1\}$ and $0^{\nabla}=1$,  and $u^{\Delta}=b, v^{\Delta}=a, a^{\Delta}=v, b^{\Delta}=u, w^{\Delta}=1, 1^{\Delta}=0$ and $0^{\Delta}=1$. By Proposition \ref{pro: irre} $(K_{7}; \vee, \wedge, ^{\Delta}, ^{\nabla}, 0, 1)$ is a WDL. 
		
		One can check that $S(K_{7})=\{0, u, v, 1\}$ and
		$\overline{S}(K_{7})=\{0, u,v, a, b, 1\}$ that are respectively a Boolean algebra and an ortholattice with $S(K_{7})\subseteq \overline{S}(K_{7})$.
		
		\item Consider the lattice $\mathcal{L}_{7}$, with $L_{7} = \{0, u, v, a, b, w, 1\}$, represented in Fig.3. One has $J(L_{7})=\{d, b, a, c\}$ and $M(L_{7})=\{e, a, b, c\}$.  
		The unary operations $^{\Delta}$ and $^{\nabla}$ are defined by the table (T:3).  
		The algebra $(L_{7}; \vee, \wedge, ^{\Delta}, ^{\nabla}, 0,1)$ is a WDL, with  
		$S(L_{7}) = \{0, 1, u, v\}$ and $\overline{S}(L_{7}) = \{0, 1, a, b\}$ that form    both Boolean algebras.
		
		\item Let $M_{7}$ and $N_{5}$ be the bounded lattices with 7 (resp. 5) elements, presented in Fig. 4 and Fig. 5. $J(M_{7})=\{d,b; a, c\}$, $M(M_{7})=\{e, a, b, c\}$.
		Using Proposition \ref{propp:wdlv} we define  the unary operations $^{\Delta}$ and $^{\nabla}$ given in (T:4) for $M_{7}$. 
		The algebra $\mathcal{M}_{7}$ is a WDL with skeletons  
		$\overline{S}(M_{7}) = \{0, 1, e, c\}$ and $S(M_{7}) = \{0, d, c, 1\}$, which are Boolean algebras.  Moreover, $\mathcal{L}_{5}$ is a WDL that is a subalgebra of $\mathcal{M}_{7}$. 
		Furthermore, $\overline{S}(M_{7}) \cap S(M_{7}) = \{0, c, 1\}$ is a three-element chain, and $S(M_{7}) \cup \overline{S}(M_{7}) = N_{5}$. 
	\end{enumerate}
\end{example}
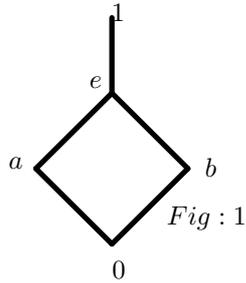
\begin{center}
	(T:1)~$\quad$	\begin{tabular}{|c|c|c|c|c|c|c|}
		\hline
		$x$ & $0$ & $u$ & $v$ & $a$ & $b$ & $1$ \\
		\hline
		$x^{\Delta}$ & $1$ & $b$ & $1$ & $b$ & $u$ &$0$ \\
		\hline
		$x^{\nabla}$ & $1$ & $v$ & $u$ & $0$ & $0$ & $0$\\
		\hline
	\end{tabular}$\quad$ 	(T: 2)$\quad$ \begin{tabular}{|c|c|c|c|c|c|c|}
		\hline
		$x$ & $0$ & $u$ & $v$ & $a$ & $b$ & $1$ \\
		\hline
		$x^{\Delta}$ & $1$ & $b$ & $1$ & $b$ & $u$ &$0$ \\
		\hline
		$x^{\nabla}$ & $1$ & $b$ & $u$ & $0$ & $u$ & $0$\\
		\hline
	\end{tabular}\\\vspace{0.3cm}  (T: 3)$\quad$\begin{tabular}{|c|c|c|c|c|c|c|c|}\hline
		$x$ & $0$ & $u$ & $v$ & $a$ & $b$ &$w$& $1$ \\
		\hline
		$x^{\Delta}$ & $1$ & $1$ & $1$ & $b$ & $a$ & $1$ &$0$ \\
		\hline
		$x^{\nabla}$ & $1$ & $v$ & $u$ & $0$ & $0$ & $0$ & $0$\\
		\hline
	\end{tabular}~(T: 4) \begin{tabular}{|c|c|c|c|c|c|c|c|}
		\hline
		$x$ & $0$ & $d$ & $a$ & $b$ & $e$ &$c$& $1$ \\
		\hline
		$x^{\Delta}$ & $1$ & $1$ & $1$ & $1$ & $c$ & $e$ &$0$ \\
		\hline
		$x^{\nabla}$ & $1$ & $c$ & $0$ & $0$ & $0$ & $d$ & $0$\\
		\hline
	\end{tabular}
\end{center}

\begin{minipage}{5cm}
	\begin{tikzpicture}[line cap=round,line join=round,>=triangle 45,x=1cm,y=1cm]
		\clip(-1.3,-1.56) rectangle (4.02,2.2);
		\draw [line width=2pt] (1,2)-- (1,1);
		\draw [line width=2pt] (1,1)-- (0,0);
		\draw [line width=2pt] (0,0)-- (1,-1);
		\draw [line width=2pt] (1,-1)-- (2,0);
		\draw [line width=2pt] (2,0)-- (1,1);
		\draw (0.88,-1.1) node[anchor=north west] {$0$};
		\draw (0.87,2.3) node[anchor=north west] {$1$};
		\draw (2.10,0.26) node[anchor=north west] {$b$};
		\draw (-0.48,0.28) node[anchor=north west] {$a$};
		\draw (0.58,1.34) node[anchor=north west] {$e$};
		\draw (1.6,-0.38) node[anchor=north west] {$Fig: 1$};
	\end{tikzpicture}
\end{minipage}
\begin{minipage}{10cm}
	
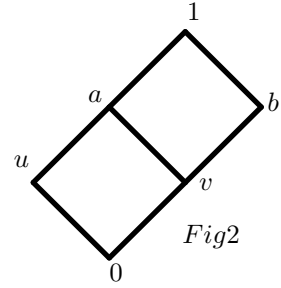
\begin{tikzpicture}[line cap=round,line join=round,>=triangle 45,x=1cm,y=1cm]
		\clip(0.3,-1.0) rectangle (10.58,3.4);
		\draw [line width=2pt] (4,1)-- (5,2);
		\draw [line width=2pt] (5,2)-- (6,3);
		\draw [line width=2pt] (6,3)-- (7,2);
		\draw [line width=2pt] (7,2)-- (6,1);
		\draw [line width=2pt] (6,1)-- (5,2);
		\draw [line width=2pt] (6,1)-- (5,0);
		\draw [line width=2pt] (5,0)-- (4,1);
		\draw (5.9,3.5) node[anchor=north west] {$1$};
		\draw (6.96,2.32) node[anchor=north west] {$b$};
		\draw (4.6,2.34) node[anchor=north west] {$a$};
		\draw (3.62,1.48) node[anchor=north west] {$u$};
		\draw (6.06,1.18) node[anchor=north west] {$v$};
		\draw (4.88,0.05) node[anchor=north west] {$0$};
		\draw (5.84,0.56) node[anchor=north west] {$Fig2$};
	\end{tikzpicture}
\end{minipage}

\begin{remark}\label{rek:soubclasse wdl}
	\begin{enumerate}
		\item Each $n$-elements chain $\mathcal{C}_{n}$ with $n \geq 3$ satisfies $\{0, 1\} = S(\mathcal{C}_{n}) = \overline{S}(\mathcal{C}_{n}) \neq \mathcal{C}_{n}$. Moreover, if a WDL  $L$ is a Boolean algebra, then $S(L) = \overline{S}(L) = L$.
		
		\item For each WDL $L \in \{\mathcal{N}_{5}, \mathcal{M}_{7}, L_{6}\}$, one has $\overline{S}(L) \neq S(L)$; both $S(L)$ and $\overline{S}(L)$ are Boolean algebras, and their intersection 
		$
		S(L) \cap \overline{S}(L) = \{0, \star, 1\}
		$ 
		is a three-element chain, which is neither an ortholattice nor a Boolean algebra.
		
		\item For the WDL $\mathcal{M}_{4,2}$, both $\overline{S}(\mathcal{M}_{4,2})$ and $S(\mathcal{M}_{4,2})$ are Boolean algebras, with $\overline{S}(\mathcal{M}_{4,2}) \subseteq S(\mathcal{M}_{4,2})$, and 
		$
		\mathcal{M}_{4,2} \neq S(\mathcal{M}_{4,2}) \cup \overline{S}(\mathcal{M}_{4,2}).
		$
		
		\item For each WDL $L \in \{L_{7}, L_{5}\}$, both $\overline{S}(L)$ and $S(L)$ are Boolean algebras, with 
		$	S(L) \cap \overline{S}(L) = \mathcal{C}_{2}.$
		Furthermore, 
		$
		\overline{S}(L_{7}) \cup S(L_{7}) \neq L_{7}, \quad \text{while} \quad N_{5} = \overline{S}(N_{5}) \cup S(N_{5}).$
		\item The WDL $\mathcal{L}_{16}$ presented in Fig.~2, page 11 of \cite{kwdd} is distributive, and its skeletons are ortholattices that are not Boolean algebras (see Fig.~4, page 22 of \cite{kwdd}). The intersection of these ortholattices forms a four-element Boolean algebra with universe 
		$
		\{0, c, a^{\nabla}, 1\}.
		$
	\end{enumerate}
\end{remark}
In the proof of Corollary. 3.3.8 in \cite{kwdd}, it is proved that $B(L)$ is  the largest subalgebra of $L$ contained in $S(L)$ and in $\overline{S}(L)$.  Notice that the inclusion $B(L)\subseteq S(L)\cap \overline{S}(L)$ may be strict, see for example the WDL $\mathcal{M}_{7}$.\par	
The previous remark allows us to introduce new subclasses within the class of WDLs based on the properties of their skeletons and dual skeletons.

\begin{minipage}{6cm}
	\begin{tikzpicture}[]
		\clip(-4.17,-3.7) rectangle (3.41,3);
		\draw [line width=2pt] (-2,1)-- (-1,2);
		\draw [line width=2pt] (-1,2)-- (0,1);
		\draw [line width=2pt] (0,1)-- (-1,0);
		\draw [line width=2pt] (-1,0)-- (-2,1);
		\draw [line width=2pt] (-2,1)-- (-2,-1);
		\draw [line width=2pt] (-2,-1)-- (-1,-2);
		\draw [line width=2pt] (-1,-2)-- (-1,0);
		\draw [line width=2pt] (0,1)-- (0,-1);
		\draw [line width=2pt] (0,-1)-- (-1,-2);
		\draw (-1.13,-2.16) node[anchor=north west] {0};
		\draw (-2.5,-0.68) node[anchor=north west] {$u$};
		\draw (0.17,-0.56) node[anchor=north west] {$v$};
		\draw (-1.3,0.64) node[anchor=north west] {$w$};
		\draw (-1.07,2.5) node[anchor=north west] {$1$};
		\draw (-2.47,1.44) node[anchor=north west] {$a$};
		\draw (0.1,1.26) node[anchor=north west] {$b$};
		\draw (-1.49,-2.62) node[anchor=north west] {$Fig ~3$};
\end{tikzpicture}\end{minipage}\begin{minipage}{10cm}
	\begin{tikzpicture}[line cap=round,line join=round,>=triangle 45,x=1cm,y=1cm]
		\clip(-1.3,-2.02) rectangle (15.02,3.7);
		\draw [line width=2pt] (-1,1)-- (0,2);
		\draw [line width=2pt] (0,2)-- (1,3);
		\draw [line width=2pt] (1,3)-- (2,1);
		\draw [line width=2pt] (0,2)-- (1,1);
		\draw [line width=2pt] (1,1)-- (0,0);
		\draw [line width=2pt] (0,0)-- (-1,1);
		\draw [line width=2pt] (0,0)-- (1,-1);
		\draw [line width=2pt] (1,-1)-- (2,1);
		\draw [line width=2pt] (5,3)-- (4,2);
		\draw [line width=2pt] (4,2)-- (4,1);
		\draw [line width=2pt] (4,1)-- (5,0);
		\draw [line width=2pt] (5,0)-- (6,1);
		\draw [line width=2pt] (6,1)-- (5,3);
		\draw (0.88,3.58) node[anchor=north west] {$1$};
		\draw (4.86,3.58) node[anchor=north west] {$1$};
		\draw (4.8,0.02) node[anchor=north west] {$0$};
		\draw (0.82,-1) node[anchor=north west] {$0$};
		\draw (2,1.28) node[anchor=north west] {$c$};
		\draw (6.0,1.36) node[anchor=north west] {$c$};
		\draw (3.64,2.36) node[anchor=north west] {$e$};
		\draw (-0.5,2.36) node[anchor=north west] {$e$};
		\draw (3.6,1.28) node[anchor=north west] {$d$};
		\draw (-0.4,0.1) node[anchor=north west] {$d$};
		\draw (1.02,1.28) node[anchor=north west] {$b$};
		\draw (-1.35,1.4) node[anchor=north west] {$a$};
		\draw (-0.88,-0.6) node[anchor=north west] {$Fig: 4$};
		\draw (5.68,0.4) node[anchor=north west] {$Fig 5$};
	\end{tikzpicture}
\end{minipage}
	\begin{minipage}{10cm}
    	\definecolor{ttqqqq}{rgb}{0.2,0,0}
	\definecolor{xdxdff}{rgb}{0.49019607843137253,0.49019607843137253,1}
	\definecolor{ududff}{rgb}{0.30196078431372547,0.30196078431372547,1}
	\begin{tikzpicture}[line cap=round,line join=round,>=triangle 45,x=1cm,y=1cm]
		\clip(-4.89,-5.78) rectangle (5.14,4.06);
		\draw [line width=2pt] (-4,1.02)-- (-3,0);
		\draw [line width=2pt] (-3,0)-- (-2.06,-1);
		\draw [line width=2pt] (-1,-2)-- (0,-1);
		\draw [line width=2pt] (0,-1)-- (1,0);
		\draw [line width=2pt] (1,0)-- (2,1);
		\draw [line width=2pt] (2,1)-- (1,2);
		\draw [line width=2pt] (1,2)-- (0,3);
		\draw [line width=2pt] (0,3)-- (-1,4);
		\draw [line width=2pt] (-1,4)-- (-2,3);
		\draw [line width=2pt] (-2,3)-- (-3.01,1.98);
		\draw [line width=2pt] (-3.01,1.98)-- (-4,1.02);
		\draw [line width=2pt] (-3,0)-- (-2.02,0.99);
		\draw [line width=2pt] (-2.02,0.99)-- (-1,2);
		\draw [line width=2pt] (-1,2)-- (0,3);
		\draw [line width=2pt] (-2.06,-1)-- (-1,0);
		\draw [line width=2pt] (-1,0)-- (0,1.06);
		\draw [line width=2pt] (0,1.06)-- (1,2);
		\draw [line width=2pt] (-3.01,1.98)-- (-2.02,0.99);
		\draw [line width=2pt] (-2.02,0.99)-- (-1,0);
		\draw [line width=2pt] (-1,0)-- (0,-1);
		\draw [line width=2pt] (-2,3)-- (-1,2);
		\draw [line width=2pt] (-1,2)-- (0,1.06);
		\draw [line width=2pt] (1,0)-- (0,1.06);
		\draw (-1.09,4.47) node[anchor=north west] {$1$};
		\draw (-2.53,3.23) node[anchor=north west] {$e$};
		\draw (-3.59,2.21) node[anchor=north west] {$d$};
		\draw (-4.5,1.16) node[anchor=north west] {$c$};
		\draw (-3.6,0.23) node[anchor=north west] {$b$};
		\draw (0.14,3.5) node[anchor=north west] {$c^{\Delta}$};
		\draw (1.19,2.49) node[anchor=north west] {$d^{\Delta}$};
		\draw (2.21,1.22) node[anchor=north west] {$e^{\Delta}=a^{\nabla}$};
		\draw (1.20,0.31) node[anchor=north west] {$b^{\nabla}$};
		\draw (-0.17,1.56) node[anchor=north west] {$v$};
		\draw (-1.15,0.57) node[anchor=north west] {$u$};
		\draw (-0.01,-0.81) node[anchor=north west] {$c^{\nabla}$};
		\draw (-1.09,-1.98) node[anchor=north west] {$0$};
		\draw [line width=2pt] (-2.06,-1)-- (-1,-2);
		\draw (-2.38,-1.14) node[anchor=north west] {$a$};
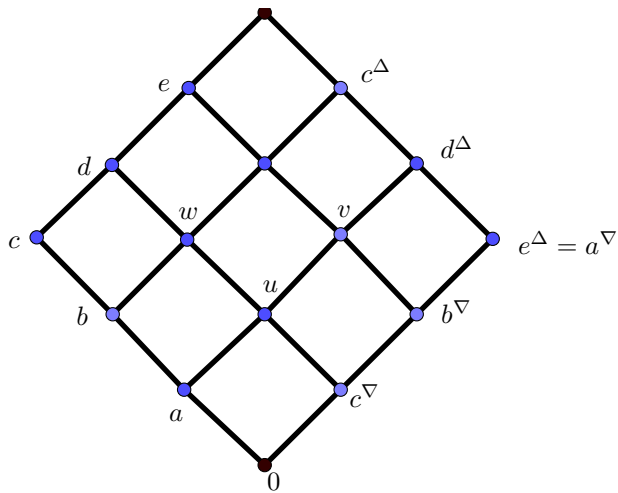
		\draw (-4.86,-2.36) node[anchor=north west] {~$Fig. 6$. \tiny A dicomplementation~on~the product~of~two~4~ element~chains.};
		\draw (-2.26,1.54) node[anchor=north west] {$w$};
		\begin{scriptsize}
			\draw [fill=ududff] (-4,1.02) circle (2.5pt);
			\draw [fill=xdxdff] (-3,0) circle (2.5pt);
			\draw [fill=ududff] (-2.06,-1) circle (2.5pt);
			\draw [fill=ttqqqq] (-1,-2) circle (2.5pt);
			\draw [fill=xdxdff] (0,-1) circle (2.5pt);
			\draw [fill=xdxdff] (1,0) circle (2.5pt);
			\draw [fill=ududff] (2,1) circle (2.5pt);
			\draw [fill=ududff] (1,2) circle (2.5pt);
			\draw [fill=xdxdff] (0,3) circle (2.5pt);
			\draw [fill=ttqqqq] (-1,4) circle (2.5pt);
			\draw [fill=ududff] (-2,3) circle (2.5pt);
			\draw [fill=ududff] (-3.01,1.98) circle (2.5pt);
			\draw [fill=ududff] (-2.02,0.99) circle (2.5pt);
			\draw [fill=ududff] (-1,2) circle (2.5pt);
			\draw [fill=ududff] (-1,0) circle (2.5pt);
			\draw [fill=xdxdff] (0,1.06) circle (2.5pt);
		\end{scriptsize}
	\end{tikzpicture}
    \end{minipage}\begin{minipage}{7cm}
	\begin{tikzpicture}[line cap=round,line join=round,>=triangle 45,x=1cm,y=1cm]
		\clip(-4.08,-1.5) rectangle (5.78,3.25);
		\draw [line width=2pt] (-1,3)-- (-2,2);
		\draw [line width=2pt] (-2,0)-- (-1,-1);
		\draw [line width=2pt] (-1,-1)-- (0,0);
		\draw [line width=2pt] (0,2)-- (-1,3);
		\draw [line width=2pt] (-1,1)-- (-2,2);
		\draw [line width=2pt] (-1,1)-- (0,2);
		\draw [line width=2pt] (0,2)-- (-1,1);
		\draw [line width=2pt] (-2,0)-- (-1,1);
		\draw [line width=2pt] (-1,1)-- (0,0);
		\draw [line width=2pt] (0,2)-- (1,1);
		\draw [line width=2pt] (1,1)-- (0,0);
		\draw [line width=2pt] (-2,2)-- (-3,1);
		\draw [line width=2pt] (-3,1)-- (-2,0);
		\draw (-0.98,3.33) node[anchor=north west] {$1$};
		\draw (-1.28,1.63) node[anchor=north west] {$w$};
		\draw (-2.48,2.17) node[anchor=north west] {$a$};
		\draw (0.11,2.30) node[anchor=north west] {$d$};
		\draw (1.1,1.37) node[anchor=north west] {$u$};
		\draw (-3.36,1.5) node[anchor=north west] {$b$};
		\draw (-2.52,0.1) node[anchor=north west] {$c$};
		\draw (0.32,0.35) node[anchor=north west] {$v$};
		\draw (-1.18,-0.94) node[anchor=north west] {$0$};
		\draw (-0.46,-0.53) node[anchor=north west] {$Fig: 7$};
	\end{tikzpicture}
\end{minipage}
\subsection{Some new classes of WDLs}~\\
In \cite{kwdd2}, Kwuida et al. provided an overview of generalizations of Boolean algebras (that is, lattices enriched with unary operations satisfying certain properties of Boolean negation). They established a relationship between known extensions of Boolean algebras and the class of WDLs, which also generalizes Boolean algebras. At the end of their work, they claimed that the  methods used in double $p$-algebras can help in understanding the structure of WDLs and suggested finding a precise description of the structure of skeletons and dual skeletons in WDLs. This suggestion motivates us to first describe various types of WDLs according to the properties verified by their skeletons and dual skeletons.  

Following Remark \ref{rek:soubclasse wdl}, we identify some new subclasses of  WDLs based on the structure of their skeletons and dual skeletons. Recall that, distributive double $p$-algebras and  Boolean algebras form the known subclasses of WDLs.

\begin{itemize}
	\item[(1)] 
	A WDL in which $(L;\vee, \wedge, 0,1)$ is a bounded distributive lattice is called a \emph{distributive WDL}. The classe of WDLs  contains many subclasses  based on the structure of the skeleton and the  dual skeleton.
	\item[(2)] A WDL $L$ is called \textbf{Boolean} if it satisfies the equalities $S(L) = \overline{S}(L) = B(L)$, for example, the WDL  $\mathcal{P}_{6}$ is Boolean.
	
	\item[(3)] A WDL $L$ is called \textbf{S-Boolean} if $S(L)$ and $\overline{S}(L)$ are Boolean algebras (even if $S(L) \neq \overline{S}(L)$), for example, $\mathcal{M}_{7}$ is S-Boolean, but not Boolean,since  $S(\mathcal{M}_{7})\neq\overline{S}(\mathcal{M}_{7})$ and the skeleton and the dual skeleton are Boolean algebras.
	
	\item[(4)] A WDL $L$ is called \textbf{pure} if $L = S(L) \cup \overline{S}(L)$. (For example, $N_{5}$ is pure, non-distributive, and S-Boolean).
	\item[(5)] A WDL is called \textbf{weak S- Boolean} if $S(L)$ or $\overline{S}(L)$ is a Boolean algebra (for example, $K_{7})$).
\end{itemize}
For a pure WDL, one can further distinguish:
(i) those that are not Boolean algebras (for example $\mathcal{L}_{5}$), 
(ii) those whose skeletons are not Boolean algebras,
(iii)  those whose skeletons are Boolean algebras. example $\mathcal{L}_{5}$.\\
For the WDLs that are not pure, one has the case where: $L$ is Boolean but not pure, for  example,  $n$-chains and trivial WDLs.\par 
In the case where $S(L)$ and $\overline{S}(L)$ are ortholattices that are not all Boolean algebras, the following cases can occur:\\
(i) $S(L) \cap \overline{S}(L) = \{0, \star, 1\}$ with $\star^{\Delta}\neq\star^{\nabla}$,\\ (ii) $S(L)\cap \overline{S}(L)=B(L)$.

\begin{example}
	\begin{enumerate}
		\item  For each $n \geq 3$, $\mathcal{C}_{n}$ as well as each Boolean algebra is a Boolean WDL. In particular,  $\mathcal{P}_{6}$ is also a Boolean WDLs.
		\item The algebras $M_{4,2}$, $\mathcal{N}_{5}$, $\mathcal{L}_{7}$, $\mathcal{M}_{7}$, $\mathcal{L}_{6}$, and $\mathcal{P}_{6}$ are $S$-Boolean WDLs.
		\item The algebras $\mathcal{M}_{4,2}$, $\mathcal{M}_{6}$, $\mathcal{M}_{7}$,  $\mathcal{L}_{7}$ and   $\mathcal{C}_{n}$, $n \geq 3$ are   $S$-Boolean WDLs that are not pure.
		\item The WDLs $\mathcal{N}_{5}$ and $\mathcal{L}_{7}$ are non-distributive $S$-Boolean WDLs.
		\item Consider the weakly dicomplementation  on  the product of two 4 elements chains presented in Fig. 6. All elements above $u$ are sent to $0$ by $^{\nabla}$. The elements $c, b$ and $a$ are each image (of their image). The operation$^{\Delta}$ is dual of $^{\nabla}$.   One can check that $S(L)=\{a, b, c, a^{\nabla}, b^{\nabla}, c^{\nabla}, 0,1\}$ and $\overline{S}(L)=\{c, d, e, c^{\Delta}, d^{\delta}, e^{\Delta}, 0,1\}$. $\overline{S}(L_{16})$ and $S(L_{16})$ are ortholattices that are not Boolean algebras. (see Fig. 1.2 pp.11  and Fig 1.4 pp.22 in \cite{kwdd}).
		\item Let $L = B \times \mathcal{C}_{n}$ be the product of WDLs, where $B$ is a Boolean algebra and $\mathcal{C}_{n}$ is a chain of $n$-elements. Then $L$ is a weakly dicomplemented lattice that is Boolean, with $B(L) = B \times \mathcal{C}_{2}$. Note that $(a, b)^{\Delta} = (a', b^{\Delta})$ and $(a, b)^{\nabla} = (a', b^{\nabla})$, where $a \in B$ and $b \in \mathcal{C}_{n}$.
		\item More generally, if we consider a WDL $L = (L; \vee, \wedge, ^{\Delta}, ^{\nabla}, 0,1)$ and a Boolean algebra $\mathcal{B} = (B; \vee, \wedge, ', 0, 1)$, then we can define
		$^{\tilde{\Delta}}$ and $^{\tilde{\nabla}}$ on $P = B \times L$ by
		\[
		(a, b)^{\tilde{\Delta}} = (a', b^{\Delta}), \quad (a, b)^{\tilde{\nabla}} = (a', b^{\nabla}).
		\]
		Then $(P; \vee, \wedge, ^{\tilde{\Delta}}, ^{\tilde{\nabla}}, \tilde{0}, \tilde{1})$
		is a WDL with $\tilde{0} = (0,0)$, $\tilde{1} = (1,1)$. Moreover, $P$ is Boolean, distributive, pure, and proper (resp. $S$-Boolean) if and only if $L$ is Boolean (resp. $S$-Boolean).
	\end{enumerate}
\end{example}

Let $L$ be a WDL and $X$ be a non-empty set. Let $L^{X}$ denotes the set of all maps from $X$ to $L$. The operations of $L$ are extended to $L^{X}$ pointwise as follows: for $f, g \in L^{X}$ and $x \in X$:

\begin{itemize}
	\item $(f \wedge g)(x) = f(x) \wedge g(x), \quad (f \vee g)(x) = f(x) \vee g(x)$,
	\item $f^{\Delta}(x) = (f(x))^{\Delta}, \quad f^{\nabla}(x) = (f(x))^{\nabla}$.
\end{itemize}  In particular, for any $a \in L$, define $\phi_{a}: X \to L$ by $\phi_{a}(x) = a$ for all $x \in X$.
The algebra $\mathcal{D} = (L^{X}; \vee, \wedge, ^{\Delta}, ^{\nabla}, \phi_{0}, \phi_{1})$ is a WDL, called the power of the WDL $L$. The corresponding order is given by 
\[
f \leq g \;\; \Leftrightarrow \;\; f(x) \leq g(x) \quad \forall x \in X.
\]
It holds that $a \leq b \;\Leftrightarrow\; \phi_{a} \leq \phi_{b}, \forall a, b\in L$.  It is easy to check that: $$\phi_{a \wedge b} = \phi_{a} \wedge \phi_{b}, \quad \phi_{a \vee b} = \phi_{a} \vee \phi_{b},
\quad  (\phi_{a})^{\Delta} = \phi_{a^{\Delta}}, \quad (\phi_{a})^{\nabla} = \phi_{a^{\nabla}}.$$
The following proposition presents some  connections between $L$ and $L^{X}$. We consider  the following particular subsets of $L$, 
$$D(L):=\{x\in L\mid x^{\nabla}=0\}~ \text{and}~ \overline{D}(L)=\{x\in L\mid x^{\Delta}=1\}.$$
One can check that $D(L)$ (resp. $(\overline{D}(L))$ is an order filter (resp. and order ideal) of $L$.
\begin{proposition}
	The following statements hold:
	\begin{enumerate}
		\item (a) $S(L^{X}) = \{f \in L^{X} \mid f(X) \subseteq S(L)\} = (S(L))^{X}$.\\
		(b) $\overline{S}(L^{X}) = \{f \in L^{X} \mid f(X) \subseteq \overline{S}(L)\} = (\overline{S}(L))^{X}$.
		\item (i) $(B(L))^{X} = B(L^{X})$, \quad (ii) $D(L^{X}) = (D(L))^{X}$, \quad (iii) $\overline{D}(L^{X}) = (\overline{D}(L))^{X}$.
		\item $L$ is Boolean (resp. $S$-Boolean ) if and only if $L^{X}$ is Boolean (resp. $S$-Boolean).
		\item $L$ is pure if and only if $L^{X}$ is pure.
		\item $L$ is distributive if and only if $L^{X}$ is distributive.
	\end{enumerate}
\end{proposition}

\begin{proof}
	(1) (a) Let $M = \{ f \in L^{X} \mid f(X) \subseteq S(L) \}$. Clearly, $M \subseteq S(L^{X})$. For the reverse inclusion, let $f \in S(L^{X})$. Then $f \in L^X$ and $f^{\nabla\nabla} = f$. Let $x \in X$. Then 
	$
	f^{\nabla\nabla}(x) = (f(x))^{\nabla\nabla} = f(x).
	$ 
	Since $f(x)^{\nabla\nabla} \in S(L)$, we deduce that $f \in M$. Therefore, $S(L^{X}) \subseteq M$ and thus $M = S(L^{X})$. A similar argument shows that 
	$\overline{S}(L^{X}) = \{ f \in L^{X} \mid f(X) \subseteq \overline{S}(L) \},
	$
	completing the proof of (1).
	
	(2) (i) Let $\phi \in (B(L))^X$. Then $\phi: X \to B(L)$, so for all $x \in X$, $\phi(x) \in B(L)$, i.e., $\phi(x) = \phi(x)^{\Delta\nabla}$. Hence, $\phi^{\Delta\nabla} = \phi$ and $\phi: X \to L$, so $\phi \in B(L^X)$. 
	
	For the reverse inclusion, let $\phi \in B(L^X)$. Then $\phi^{\Delta\nabla} = \phi$. For any $x \in X$, we have 
	$
	\phi(x) = (\phi(x))^{\Delta\nabla},
	$
	so $\phi(x) \in B(L)$. Therefore, $\phi \in B(L)^X$. This proves the desired equality. The statements in (ii) and (iii) can be proved similarly.
	
	(3) Assume that $L$ is a Boolean WDL,  i.e. $S(L) = \overline{S}(L)$ is a Boolean algebra. We show that $S(L^{X}) = \overline{S}(L^{X})$. Since $S(L) = \overline{S}(L)$, part (1) implies the result.  
	
	For the converse, assume that $L^{X}$ is a Boolean WDL, i.e., $S(L^{X}) = \overline{S}(L^{X})$. We show that $\overline{S}(L) \subseteq S(L)$. Let $a \in \overline{S}(L)$. Then $a^{\Delta\Delta} = a$, so $\phi_a^{\Delta\Delta} = \phi_a$. It follows that $\phi_a \in \overline{S}(L^{X}) = S(L^{X})$, hence $(\phi_a)^{\nabla\nabla} = \phi_a$, i.e., $\phi_{a^{\nabla\nabla}} = \phi_a$. We deduce that $a^{\nabla\nabla} = a$, so $a \in S(L)$. Thus, $\overline{S}(L) \subseteq S(L)$. A dual argument shows that $S(L) \subseteq \overline{S}(L)$, hence $S(L) = \overline{S}(L)$. This completes the proof of (3).
	
	(4) and (5) are straightforward to verify.
\end{proof}
\begin{example}
	Consider the  WDL $\mathcal{C}_{3}$ with $C_{3}=\{0,e, 1\}$. Let $X=\{0,1\}$ and  $L_{9}=C_{3}^{X}$ be the power of $\mathcal{C}_{3}$ presented in Fig. 7. We  have   $L_{9}=\{0, a,b, c, d, u, v, w, 1\}$, with   $a=(e,1), b=(0,1), c=(0,e), d=(1,e), u=(1, 0), v=(e, 0), w=(e,e)$ and $1=(1,1)$. In addition $L_{9}$ and $\mathcal{C}_{3}$ are Boolean WDL. In fact, one has $S(L_{9})=\overline{S}(L_{9})=\{0,1, u, b\}$ which is a Boolean algebra.
\end{example}
\begin{minipage}{8cm}
	\begin{tabular}{|c|c|c|c|c|c|c|c|c|c|}
		\hline
		$x$ & $0$ & $a$ & $b$ & $c$ & $d$ &$u$& $v$ & $w$& $1$\\
		\hline
		$x^{\Delta}$ & $1$ & $u$ & $u$ & $1$ & $b$ & $b$ &$1$&$1$&$0$ \\
		\hline
		$x^{\nabla}$ & $1$ & $0$ & $u$ & $u$ & $0$ & $b$ & $b$&$0$&$0$\\
		\hline
	\end{tabular} 
\end{minipage}

	\section{Normal filters in weakly dicomplemented lattices}

In this section, we  study the class of normal filters (resp. ideals)  in arbitrary WDLs and investigate some of their algebraic properties with the aim of addressing the problem of congruences in WDLs. We show that the lattice of normal filters of a WDL  
$L$ is isomorphic to the lattice of normal ideals of 
$L$, but it is not a sublattice of $F(L)$.
Moreover, we characterize normal filters generated by a subset, and describe the joins and the meet of two normal filters (resp. ideals). Principal normal filters are also examined and their properties discussed.

\begin{definition} 
	\begin{enumerate}
		\item A filter $F$ is called a \emph{normal filter} if, for every $x \in L$, $x \in F$ implies $x^{\Delta\nabla} \in F$.
		\item An ideal $I$ is called a \emph{normal ideal} if, for every $x \in L$, $x \in I$ implies $x^{\nabla\Delta} \in I$.
	\end{enumerate}
\end{definition}	

By using (9) of Proposition~\ref{prop: axiom}, for any $x \in L$, one has
\[
(\ast\ast)\quad x^{\Delta\nabla} \leq x^{\Delta\Delta} \leq x \leq x^{\nabla\nabla} \leq x^{\nabla\Delta}.
\]
It follows that for any normal filter (resp. ideal) of $L$, $x \in F$ (resp. $x \in I$) implies that $x^{\Delta\Delta}, x^{\nabla\nabla} \in F$ (resp. $x^{\Delta\Delta}, x^{\nabla\nabla} \in I$). For the following, we consider the unary operations $\eta_{\Delta}, \eta_{\nabla}: L\to L$ and the binary operations $\overline{\wedge}, \tilde{\vee}: L^{2}\to L$ defined by: 
$$ \eta_{\Delta}(x) = x^{\Delta\nabla},~ \eta_{\nabla}(x) = x^{\nabla\Delta},~
 x \overline{\wedge} y = (x \wedge y)^{\Delta\nabla},~
 x \tilde{\vee} y = (x \vee y)^{\nabla\Delta}.$$

Given $a \in L$, we define inductively $a^{n(\Delta\nabla)}$ and $a^{n(\nabla\Delta)}$ as follows:
\begin{itemize}
	\item $a^{1(\Delta\nabla)} = a^{\Delta\nabla}$, and $a^{(k+1)(\Delta\nabla)} = \big(a^{k(\Delta\nabla)}\big)^{\Delta\nabla}$ for every $k \geq 1$.
	\item $a^{1(\nabla\Delta)} = a^{\nabla\Delta}$, and $a^{(k+1)(\nabla\Delta)} = \big(a^{k(\nabla\Delta)}\big)^{\nabla\Delta}$ for every $k \geq 1$.
\end{itemize}	
The chain	$
a \geq a^{\Delta\nabla} \geq \cdots \geq a^{n(\Delta\nabla)} \geq \cdots
$
is called the \emph{normal chain} determined by $a$.  $a\leq a^{\nabla\Delta}\leq\ldots\leq a^{n\nabla\Delta}...$ is called \emph{dual normal chain} determined by $a$. The normal chain defined by $a$ is said to be \emph{finite} if there exists $n \geq 1$ such that 
$a^{(n+1)(\Delta\nabla)} = a^{n(\Delta\nabla)}.
$

We list some properties of the maps $\overline{\wedge}$, $\tilde{\vee}$, $\eta_{\Delta}$, and $\eta_{\nabla}$ in the following remark.
\begin{remark}\label{rk:eng}
	Let $x, y \in L$. The following statements hold:
	\begin{enumerate}
		\item[(1)] 
		(i) $\eta_{\Delta}$ and $\eta_{\nabla}$ are increasing maps.
		(ii) $x \leq \eta_{\nabla}(x)$.
		(iii) $\eta_{\Delta}(x) \leq x$.				
		\item[(2)] 
		(i) $\eta_{\Delta}(x \,\overline{\wedge}\, y) = \eta_{\Delta}(x) \,\overline{\wedge}\, \eta_{\Delta}(y)$.
		(ii) $\eta_{\Delta}(x \wedge y) = \eta_{\Delta}(x) \wedge \eta_{\Delta}(y)$.
		
		\item[(3)] 
		(i) $\eta_{\nabla}(x \vee y) = \eta_{\nabla}(x) \vee \eta_{\nabla}(y)$.
		(ii) $\eta_{\nabla}(x \,\tilde{\vee}\, y) = \eta_{\nabla}(x) \,\tilde{\vee}\, \eta_{\nabla}(y)$.
		\item[(4)] $L$ is a Boolean algebra if and only if $\eta_{\Delta} \circ \eta_{\Delta} = \eta_{\Delta}$ and $\eta_{\nabla} \circ \eta_{\nabla} = \eta_{\nabla}$.
		
		\item[(5)] $x \,\overline{\wedge}\, y \leq x \,\overline{\sqcap}\, y \leq x \wedge y$, and $x \vee y \leq x \sqcup y \leq x \,\tilde{\vee}\, y$.
		
		\item[(6)] $x^{\Delta} \,\overline{\wedge}\, x = 0$ and $x \,\tilde{\vee}\, x^{\nabla} = 1$.
		
		\item[(7)] If $x \leq y$ and $a \leq b$, then $a \,\gamma\, x \leq b \,\gamma\, y$ for all $\gamma \in \{\overline{\wedge}, \overline{\vee}, \overline{\sqcap}, \underline{\sqcup}, \sqcup, \sqcap\}$.
		
		\item[(8)] If $n \leq m$, then $a^{m(\Delta\nabla)} \leq a^{n(\Delta\nabla)}$.
		
		\item[(9)] $(a \wedge b)^{n(\Delta\nabla)} = a^{n(\Delta\nabla)} \wedge b^{n(\Delta\nabla)}$, and 
		$(a \vee b)^{n(\nabla\Delta)} = a^{n(\nabla\Delta)} \vee b^{n(\nabla\Delta)}$.
	\end{enumerate}
\end{remark}

\begin{proof}
	Easy to show, using definitions and (9) of Proposition \ref{prop: axiom}.
\end{proof}
One can observe that, in the  WDL $\mathcal{P}_{6}$, $F_{1}=\{1, a, u\}$ and $F_{2}=\{b, 1\}$ are normal filters and $F=F_{1}\vee F_{2}=\{v, u, a, b, 1\}$ is a filter of $L$ that is not normal. In fact, one has $v\in F$ and $v^{\Delta\nabla}=0\not\in F$, this observation allows us to introduce the join operation in the lattice of normal filters (ideals) and characterize the normal filter generated by a subset. For $a\in L$ and $X\subseteq L$, one denotes by $N[a)$ (resp. $N(a]$) the  principal normal filter (resp. normal ideal) of $L$ generated by $a$, by $N[X)$ (resp. $N(X]$) the normal filter (resp. normal ideal) of $L$ generated by $X$, and by $F_{1}\overline{\vee}F_{2}$ (resp. $I_{1}\tilde{\vee} I_{2}$)   the join in the lattice $NF(L)$ (resp.$NI(L)$) of normal filters (resp. ideals) of $L$. $NF(L)$ (resp. $NI(L)$ denotes the set of normal filters (resp. normal ideals) of $L$. 	Let $a, b\in L, X\subseteq L, F, G\in NF(L), I, J\in NI(L)$.

\begin{proposition}\label{nor:filter}
	The following statements hold:
	\begin{enumerate}
		\item $N[X) = \{x \in L \mid \exists\, x_{1}, \ldots, x_{n} \in X,\ n, m \geq 1,\ (x_{1} \wedge \ldots \wedge x_{n})^{m(\Delta\nabla)} \leq x\}$.  
		In particular,  
		$$F \overline{\vee} G = \{x \in L \mid \exists\, a \in F,\, \exists\, b \in G,\ a \overline{\wedge} b \leq x\}.$$
		
		\item $N(X] = \{x \in L \mid \exists\, x_{1}, \ldots, x_{n} \in X,\ n, m \geq 1,\ x \leq (x_{1} \vee \ldots \vee x_{n})^{m(\nabla\Delta)}\}$.  
		In particular,  
		$$I \tilde{\vee} J = \{x \in L \mid \exists\, a \in I,\, \exists\, b \in J,\ x \leq a \tilde{\vee} b\}.$$
		\item (i) $N[a) = \{x \in L \mid \exists\, n \geq 1,\ a^{n\Delta\nabla} \leq x\}$.\\  
		(ii) $N(a] = \{x \in L \mid \exists\, n \geq 1,\ x \leq a^{n\nabla\Delta}\}$.
		
		\item If $a \leq b$, then $N[b) \subseteq N[a)$ and $N(a] \subseteq N(b]$.
		\item (i) $N[a) = N[a^{\Delta\nabla})$.  
		(ii) $N(a] = N(a^{\nabla\Delta}]$.
		\item (i) $N[a) \overline{\vee} N[b) = N[a \overline{\wedge} b)$.   (ii)   $N[a)\cap N[b)\subseteq N[a\vee b)$.
		\item  (i) $N(a]\overline{\vee}N(b]=N(a\overline{\vee} b]$.  (ii)
		$N(a]\cap N(b] \subseteq N[a \overline{\vee} b)$.
		\item $N[a) \overline{\vee} N[a^{\Delta}) = L$ and $N(a] \tilde{\vee} N(a^{\nabla}] = L$.

		\item   
		(i) $N[a) = L$ iff $\exists\, n \geq 1$,  $a^{n\Delta\nabla} = 0$. (ii) $N(a]=L$ iff $\exists n\geq 1, a^{n\nabla\Delta}=1$.
	\end{enumerate}
\end{proposition}

\begin{proof}
	(1) If $X = \emptyset$, then $N[X) = \{1\}$.  
	Assume that $X \neq \emptyset$. Let  
	$$
	M = \{x \in L \mid \exists\, x_{1}, \ldots, x_{n} \in X,\ n, m \geq 1,\ (x_{1} \wedge \ldots \wedge x_{n})^{m(\Delta\nabla)} \leq x\}.
	$$
	Clearly, $X \subseteq M$, since $x^{\Delta\nabla} \leq x$ for all $x \in L$.  
	It is also clear that if $x \leq y$ and $x \in M$, then $y \in M$, by transitivity of $\leq$ and the definition of $M$.
	
	Let $x, y \in M$. Then there exist $n, m, k, l \geq 1$ and $a_{1}, \ldots, a_{n}, b_{1}, \ldots, b_{m} \in X$ such that  			
	$$(\star_{1})\quad\quad
	(a_{1} \wedge \ldots \wedge a_{n})^{k(\Delta\nabla)} \leq x 
	~	\text{and}~(\star_2)~
	(b_{1} \wedge \ldots \wedge b_{m})^{l(\Delta\nabla)} \leq y.$$
	Note that  
	$(\star_3)~
	a^{n(\Delta\nabla)} \wedge b^{n(\Delta\nabla)} = (a \wedge b)^{n(\Delta\nabla)}$ and if $n \leq m$, then  
	$ (\star_{4})\quad 
	a^{n(\Delta\nabla)} \geq a^{m(\Delta\nabla)}. (\star_{4})$ by (1) of Remark\ref{rk:eng}.	Without loss of generality, assume that $k \leq l$. Then set
	\[
	c = (a_{1} \wedge \ldots \wedge a_{n})^{l(\Delta\nabla)} \wedge (b_{1} \wedge \ldots \wedge b_{m})^{l(\Delta\nabla)}.
	\]
	We have
	$
	c \leq (a_{1} \wedge \ldots \wedge a_{n})^{k(\Delta\nabla)} \wedge (b_{1} \wedge \ldots \wedge b_{m})^{l(\Delta\nabla)} \leq x \wedge y.$
	By ($\star_3$), we get
	\[
	c = (a_{1} \wedge \ldots \wedge a_{n} \wedge b_{1} \wedge \ldots \wedge b_{m})^{l(\Delta\nabla)} \leq x \wedge y.
	\]
	Thus $x \wedge y \in M$, so $M$ is a filter of $L$.
	
	Next, we show that $M$ is a normal filter.  
	Since $\eta_{\Delta}$ is an increasing function, we have
	\[
	[(a_{1} \wedge \ldots \wedge a_{n})^{k(\Delta\nabla)}]^{\Delta\nabla}
	= (a_{1} \wedge \ldots \wedge a_{n})^{(k+1)(\Delta\nabla)} \leq x^{\Delta\nabla},
	\]
	which implies that $x^{\Delta\nabla} \in M$. Hence $M$ is a normal filter of $L$ containing $X$.
	
	Let $H$ be any normal filter of $L$ containing $X$. We show that $M \subseteq H$.  
	Let $z \in M$. Then there exist $a_{1}, \ldots, a_{n} \in X$ and $k, m \geq 1$ such that  
	\[
	(a_{1} \wedge \ldots \wedge a_{n})^{k(\Delta\nabla)} \leq z.
	\]
	Since $H$ is a normal filter and $a_{1}, \ldots, a_{n} \in X \subseteq H$, we have $z \in H$.  
	Thus $M \subseteq H$, and therefore $N[X) = M$. Now, let
	\[
	H = \{x \in L \mid \exists\, u \in F,\, \exists\, v \in G,\ u \overline{\wedge} v \leq x\}.
	\]
	A similar argument as above shows that the desired equality holds. \\The item (2) is dual to item (1).\\	
	\vspace{0.2cm}
	(3) (i) Let $K = \{x \in L \mid \exists\, n \geq 1,\ a^{n\Delta\nabla} \leq x\}$.  
	Clearly, $a, a^{\Delta\nabla} \in K$. Moreover, if $x \leq y$ and $x \in K$, then $y \in K$.  
	Assume that $x, y \in K$. Then there exist $n, m \geq 1$ such that $a^{n\Delta\nabla} \leq x$ and $a^{m\Delta\nabla} \leq y$.  
	If $n \leq m$, then $a^{m\Delta\nabla} \leq x \wedge y$, hence $x \wedge y \in K$.  
	Clearly, $K$ is the least normal filter containing $a$, and thus $K = N[a)$. (ii) can be proved by a dual argument.\\
	(4) and (5) are easy to verify.\\
	(6)(i) We now show that $N[a) \overline{\vee} N[b) = N[a \wedge b)$.  
	By (4), we have $N[a) \overline{\vee}N[b)\subseteq N[a\wedge b)$. For the reverse inclusion. Let $z\in N[a\wedge b)$, then there exists $n\geq 1$ such that $(a\wedge b)^{n\Delta\nabla}\leq z$, and by (9) of Remark \ref{rk:eng} one has $a^{n\Delta\nabla}\overline{\wedge}b^{n\Delta\nabla}\leq a^{n\Delta\nabla}\wedge b^{n\Delta\nabla}=(a\wedge b)^{n\Delta\nabla}\leq z$. Using (2) one deduces that $z\in N[a)\overline{\vee} N[b)$, so $N[a\wedge b)\subseteq N[a)\overline{\vee} N[b)$.  Hence 6. (i) holds by 5.\\
	Since $a \vee b \in N[a) \cap N[b)$, by the minimality of $N[a \vee b)$ as the normal filter containing $a \vee b$,  
	we have $N[a \vee b) \subseteq N[a) \cap N[b)$. \\
	Item (7) is dual to item (6). 
	Item (8)(i) holds by taking $b=a^{\Delta}$ in (6)(i) and using (6) of Remark \ref{rk:eng}. 
	Similarly one can prove (8)(ii).
	\\ For (9)(ii) using (5)(i) and (3)(i) we are done. (9) (i) can be proved similarly.

\end{proof}

Given a filter $F$ and an ideal $J$ of $L$, Kwuida in \cite{kwdd} considers the sets  
\[
J(F) := \{z \in L \mid z \leq x^{\nabla} \text{ for some } x \in F\} 
\quad \text{and} \quad 
F(J) := \{z \in L \mid z \geq x^{\Delta} \text{ for some } x \in J\}.
\]
We have the following known results.

\begin{proposition}\label{prop knormal}
	(\cite{kwdd}, Props. 2.4.3 and 2.4.4)  
	Let $x \in L$, $F \in NF(L)$, and $J \in NI(L)$.  
	The following statements hold:
	\begin{enumerate}
		\item[(i)] $J(F)$ is a normal ideal of $L$. Moreover, 
		if $x \in F$, then $x^{\Delta}, x^{\nabla} \in J(F)$.				
		\item[(i)] $F(J)$ is a normal filter of $L$.  
		Moreover, if $x \in J$, then $x^{\nabla}, x^{\Delta} \in F(J)$.
		
		\item $F = F(J(F))$ and $J = J(F(J))$.
	\end{enumerate}
\end{proposition}	
In the next, we 
show that the lattice of normal filter of $L$ is isomorphic to the lattice of normal ideals of $L$. To that aim, we consider the following two  maps:  
\[
\Psi: NI(L) \to NF(L), \quad J \mapsto \Psi(J) = F(J), 
~ \text{and}~ 
\Phi: NF(L) \to NI(L),~ F \mapsto \Phi(F) = J(F).
\]
Clearly, by items (1) and (2) of Proposition~\ref{prop knormal}, the maps $\Psi$ and $\Phi$ are well defined.  
We now state the following useful lemma.

\begin{lemma}\label{lem iso nfilter}
	Let $F, G \in NF(L)$ and $I, J \in NI(L)$. The following statements hold:
	\begin{enumerate}
		\item If $I \subseteq J$, then $\Psi(I) \subseteq \Psi(J)$.
		
		\item (i) $\Psi(I \cap J) = \Psi(I) \cap \Psi(J)$.  
		(ii) $\Psi(I \tilde{\vee} J) = \Psi(I) \overline{\vee} \Psi(J)$.
		
		\item $\Psi$ is a bijective map, and its inverse is $\Phi$.
	\end{enumerate}
\end{lemma}
\begin{proof}
	Item (1) is straightforward to verify.
	\smallskip
	\noindent
	(2)  
	For (i), it is clear that $\Psi(I \cap J) \subseteq \Psi(I) \cap \Psi(J)$ by using (1).  	Let $z \in \Psi(I) \cap \Psi(J)$. Then there exist $u \in I$ and $v \in J$ such that  
	$z \geq u^{\Delta}$ and $z \geq v^{\Delta}$. It follows that  
	$
	z \geq u^{\Delta} \vee v^{\Delta} = (u \wedge v)^{\Delta}
	$
	by (4) of Proposition~\ref{prop: axiom}.  
	Since $u \wedge v \in I \cap J$, we deduce that $z \in \Psi(I \cap J)$.  
	Hence, $\Psi(I) \cap \Psi(J) \subseteq \Psi(I \cap J)$. Therefore $\Psi(I\cap J)=\Psi(I)\cap\Psi(J)$.
	
	\smallskip
	\noindent
	We have $\Psi(I), \Psi(J) \subseteq \Psi(I \tilde{\vee} J)$, 
	so $\Psi(I) \overline{\vee} \Psi(J) \subseteq \Psi(I \tilde{\vee} J)$.
	
	For the reverse inclusion, let $z \in \Psi(I \tilde{\vee} J)$.  
	Then there exists $x \in I \tilde{\vee} J$ such that
	$z \geq x^{\Delta} \quad(\star)$.
	Since $x \in I \tilde{\vee} J$, there exist $u \in I$ and $v \in J$ such that
	$x \leq (u \vee v)^{\nabla\Delta} \quad (\star\star)$.
	By (4') of Proposition~\ref{prop: axiom}, from $(\star\star)$ we obtain
	\[
	x^{\Delta} \geq (u \vee v)^{\nabla\Delta\Delta} \geq (u \vee v)^{\nabla\Delta\nabla}
	= (u^{\nabla} \wedge v^{\nabla})^{\Delta\nabla},
	\]
	because  $x^{\nabla} \leq x^{\Delta}$ always holds in $L$.
	Since  $I$ and $J$ are normal ideals, we have $u^{\nabla\Delta} \in I$ and $v^{\nabla\Delta} \in J$.  
	Moreover, $u^{\nabla} \geq u^{(\nabla\Delta)\Delta}$ and $v^{\nabla} \geq v^{\nabla\Delta\Delta}$.  
	Thus,
	$
	z \geq x^{\Delta} \geq (u^{\nabla} \wedge v^{\nabla})^{\Delta\nabla},
	$
	and since $u^{\nabla} \in \Psi(I)$ and $v^{\nabla} \in \Psi(J)$, we deduce that  
	$z \in \Psi(I) \overline{\vee} \Psi(J)$.  
	Hence, $\Psi(I \tilde{\vee} J) = \Psi(I) \overline{\vee} \Psi(J)$, and (ii) holds.
	
	\smallskip
	\noindent
	(3) We now show that $\Psi$ is a bijection.
	
	\textit{Injectivity.}  
	Assume that $\Psi(I) = \Psi(J)$.
	For the inclusion $I \subseteq J$, let $u \in I$.  
	Then $u^{\nabla} \in \Psi(I) = \Psi(J)$ by (2) of Proposition~\ref{prop knormal}.  
	Hence, there exists $v \in J$ such that $u^{\nabla} \geq v^{\Delta} \geq v^{\nabla}$.  
	It follows that $u \leq u^{\nabla\nabla} \leq v^{\nabla\nabla}$.  
	Since $J$ is a normal ideal, $v^{\nabla\nabla} \in J$.  
	Moreover, as $J$ is a downset, we conclude that $u \in J$. 
	Therefore, $I \subseteq J$.  
	A similar argument shows that $J \subseteq I$, hence $I = J$.  
	Thus, $\Psi$ is injective.
	
	\textit{Surjectivity.}  
	Let $F$ be a normal filter of $L$. We seek a normal ideal $I$ of $L$ such that $\Psi(I) = F$.  
	Let $I = J(F)$.  
	By (1) of Proposition~\ref{prop knormal}, $I$ is a normal ideal of $L$ and by (3) of Proposition~\ref{prop knormal}, $F=\Psi(I)$.
	Therefore, $\Psi$ is surjective.  
	Finally, by (3) of Proposition~\ref{prop knormal}, $\Psi$ and $\Phi$ are inverses of each other.
\end{proof} 

The next result shows that normal ideals and normal filters occur in pairs in WDLs.

\begin{theorem}\label{theo:isonor}
	The lattice of normal ideals of $L$ is isomorphic to the lattice of normal filters of $L$, via the map $\Psi$, with inverse given by the map $\Phi$.
\end{theorem}

\begin{proof}
	By (1), (2), and (3) of Lemma \ref{lem iso nfilter}, the result follows.
\end{proof}

\begin{proposition}\label{prop:characnf}
	(Characterization of normal filters in WDLs)
	Let $F$ be an order filter of $L$.  The following statments are equivalent.
	
	\begin{enumerate}
		
		\item For any $x, y\in L$, $x, y\in F$ implies $x\overline{\wedge} y\in F$.
		\item $F$ is a normal filter.
	\end{enumerate}
\end{proposition}
\begin{proof}
	(1)$\Rightarrow$ (2) Assume that (1) holds.  Let $x, y\in F$, then $x\overline{\wedge} y\leq x\wedge y$, as $F$is an order filter, we have $x\wedge y\in F$. In addition, $x^{\Delta\nabla}=x\overline{\wedge} x\in F$ by assumption, hence $F$ is a normal filter.\\
	(2) $\Rightarrow$ (1) Assume that $F$ is a normal filter. Let $x, y\in F$, then $x\wedge y\in F$, as $F$ is normal, we have $ (x\wedge y)^{\Delta\nabla}=x\overline{\wedge }y\in F$ and we are done.	
\end{proof}

Now let $M$ be any  subalgebra of $L$ and consider  $G\in NF(M)$. We denote by $F_{G}$ the subset of $L$ defined by   $$F_{G}:=\{x\in L\mid \exists u\in G, u\leq x\}.$$
Recall that $B(L)$  is a subalgebra of $L$  verifying $x^{\Delta}=x^{\nabla}$, for any $x\in B$,and a subset $G$ of $L$ is a base of some filter $F$ if $F=\{x\in L\mid \exists v\in G, v\leq x\}$. 
\begin{lemma}\label{lemmiso}
	If  $M$ is  a subalgebra of $L$ and  $G, H\in NF(M)$. Then: 
	\begin{enumerate}
		\item (i) $F_{G}$ is a normal filter of $L$. (ii)
		$G\subseteq H \Leftrightarrow F_{G}\subseteq F_{H}$.
		\item If $F\in NF(L)$, then $F\cap M=E\in NF(M)$ and $F_{E}\subseteq F$. Moreover, if $M=B(L)$ and for each $f\in F$, there exists $n\geq 1$ such that $f^{n\Delta\nabla}=f^{(n+1)\Delta\nabla}$,  then $F=F_{E}$.
		\item Each filter of $B(L)$ is a base of some normal filter of $L$.
		\item $F_{G}=[G)_{L}$.
		$F_{G}\cap F_{H}=F_{G\cap H}$, $F_{G}\overline{\vee} F_{H}=F_{G\overline{\vee} H}$.
		\item $B(M)\subseteq B(L)$.
	\end{enumerate}
\end{lemma}
\begin{proof}
	1. (i) Assume that $G$ is a normal filter in $M$. Clearly $1\in F_{G}$. Assume that $x\leq y$ and $x\in F_{G}$, then there exists $u\in G$ such that $u\leq x\leq y$, then $y\in F_{G}$ (by definition of $F_{G}$ and transitivity of $\leq$). Assume that $x, y\in F_{G}$, then there exist $u, v\in G$ such that $u\leq x$ and  $v\leq y$, since $G$ is normal and $u, v\in G$, we have $u\overline{\wedge}v\in G$ and $u\overline{\wedge} v\leq x\overline{\wedge} y$, as $u\overline{\wedge} v\in G$, we deduce that $x\overline{\wedge} y\in F_{G}$, by Proposition \ref{prop:characnf} $F_{G}$ is a normal filter of $L$.\vspace{0.2cm}\\	
	(ii) The direct implication is obvious. Assume that $F_{G}\subseteq F_{H}$. Let $x\in G$, then $x\in F_{G}=F_{H}$, then there exists $v\in H$ such that $v\leq x$, since $v\in H$ and $H$ is an order filter we deduce that $x\in H$. Hence $G\subseteq H$. Therefore the equivalence holds.\\
	2. The first part is easy to check, since $M$ is a subalgebra of $L$. Note that the result still holds for $M=B(L)$. Clearly $F_{E}\subseteq F$. Assume that $M=B(L)$ and for each $f\in F$, there exists $n\geq 1$ such that $f^{n\Delta\nabla}=f^{(n+1)\Delta\nabla}$.  The last equality shows that for  each $f\in F$, $f^{\Delta\nabla}\in B(L)\cap F=E$, then $f\in F_{E}$, and the desired equality holds.\\
	3. Follows from 1. in case where $M=B(L)$.\\
	4. Recall that$[G)_{L}$ denotes the filter of $L$ generated by $G$. (i)  By definition of $F_{G}$ one has $F_{G}\subseteq [G)_{L}$. For the reverse inclusion, let $x\in [G)_{L}$, then there exists $u_{1}, \ldots, u_{n}\in G, n\geq 1$ such that $u_{1}\wedge\ldots\wedge u_{n}\leq x$. Since $^{\Delta\nabla}$ is an incresing map, one deduces that $(u_{1}\wedge\ldots\wedge)^{\Delta\nabla}\leq x^{\Delta\nabla}\leq x$, as $G$ is a normal filter of $M$,  $(u_{1}\wedge\ldots\wedge u_{n})^{\Delta\nabla}\in G$; therefore $x\in F_{G}$ and $[G)_{L}\subseteq F_{G}$. Hence $F_{G}=[G)_{L}$.\\
	 (ii) Clearly one has $F_{G\cap H}\subseteq F_{G}\cap F_{H}$. Let $z\in F_{G}\cap F_{H}$, then there exists $u\in G, v\in H$ such that $u, v\leq z$. One has $(u\vee v)^{\Delta\nabla}\leq z^{\Delta\nabla}\leq z$ and $(u\vee v)^{\Delta\nabla}\in G\cap H$. It follows that $z\in F_{G\cap H}$. hance $F_{G}\cap F_{H}=F_{G\cap H}$. (iii) Clearly one has $F_{G}\overline{\vee} F_{H}\subseteq F_{G\overline{\vee} G}$. For the reverse inclusion, let $z\in F_{G\overline{\vee}H}$, then there exist $u\in G, v\in H$ such that $u^{\Delta\nabla}\wedge v^{\Delta\nabla}=(u\wedge v)^{\Delta\nabla}\subseteq z$. As $u\in F_{G}$ and $v\in F_{H}$, one deduces that $z\in F_{G}\overline{\vee} F_{H}$. Hence $F_{G}\overline{\vee} F_{H}=F_{G\overline{\vee}H}$.
	\\The item (5) is obvious.
	\end{proof}

The following theorem establishes the isomorphism between the lattice  $NF(L)$ and the lattice  $F(B(L))$ of filters of the Boolean center of $L$, under certain cndition. Here, $F(B(L))$ denotes the set of filters of the Boolean center $B(L)$ of $L$.
\begin{theorem}\label{iso: normal}
	\begin{enumerate}
		\item If $M$ is a subalgebra of $L$, then the  lattice of normal filters of $M$ embedds in the lattice of normal filters of $\mathcal{L}$.
		\item 
	If $L$ is a WDL in which every normal chain is finite, then $N(F(L))$ and $F(B(L))$ are isomorphic lattices.
	
	\end{enumerate}
\end{theorem}
\begin{proof}
	Let $\eta: NF(M)\to F(L), G\mapsto \eta(F)=F_{G}$.
By	1. and 4. of Lemma \ref{lemmiso} $\eta$ is an embedding of lattices.\par 
2. Since $B(L)$ is a subalgebra of $L$, by 1. $\eta: NF(B(L))\to NF(L), G\mapsto \eta(G)=F_{G}$ is an embedding.  The second part of 2. in Lemma \ref{lemmiso} shows that $\eta$ is  surjetive under the assumptions.	Hence $\eta$ is an isomorphism between $NF(L)$ and $NF(B(L))$ when the assumption holds.

\end{proof}
\begin{corollary}
	In each finite WDL $L$, $NF(L)$ and $NI(L)$ are isomorphic lattices.
\end{corollary}
\begin{proof}
	It is obvious, since, in every finite WDL, each normal chain is finite.
\end{proof}
As example, for the WDL $\mathcal{P}_{6}$, one has $S(\mathcal{P}_{6})=\overline{S}(\mathcal{P}_{6})=\{0, u, b, 1\}$ which forms a Boolean algebra. We have $NF(\mathcal{P}_{6})\cong F(B(\mathcal{P}_{6}))$.

\section{Congruences in distributive  WDLs}
The characterization of congruences provides useful tools for simplifying the study of WDLs. The description of subdirectly irreducible and simple WDLs depends on the lattice of congruences of a given WDL. Congruences generated by filters offer an effective approach to obtaining the desired characterizations. As we will see, normal filters play an important role in these characterizations.

Recall that concept algebras are examples of WDLs, although WDLs are not necessarily concept algebras. In \cite{kwdd} (Sect. 2.3.1), Kwuida initiated the study of congruences in distributive concept algebras using compatible subcontexts. 

In\cite{Tenk4}, we initiate investigation of congruences generated by S-filters in  distributive weakly complemented lattices.  Recall that a weakly complemented lattice $(L; \vee, \wedge, ^{\Delta}, 0, 1)$ is an algebra reduct of a weakly dicomplemented latiice $(L; \vee, \wedge, ^{\Delta}, ^{\nabla}, 0,1)$.

In the present work, we characterize congruences in distributive WDLs using normal filters. The resulting characterizations are more general and simpler than the existing ones for distributive concept algebras.	
\subsection{Some preliminaries results on congruences in WDLs}~\par
It is well known (see \cite{gratz2}) that an equivalence relation $\theta$ on a lattice $(L; \vee, \wedge)$ is a congruence relation on $L$ if, for all $x, y, z \in L$,	
$$(\star)\quad (x, y)\in\theta\Rightarrow (x\wedge z, y\wedge z), (x\vee z, y\vee z)\in\theta.$$
If $F \in NF(L)$, then $\Theta(F)$ denotes the smallest congruence on $L$ collapsing $F$, and $\theta(a, b)$ denotes the principal congruence relation on $L$ collapsing the pair $a, b \in L$.
For the following, we consider the binary relation $\Phi$ on $L$ defined by
$$\Phi:=\{(x, y)\in L^{2}\mid x^{\nabla}=y^{\nabla}~\text{and}~x^{\Delta}=y^{\Delta}\}.$$

It is easy to verify that $\Phi \in \mathrm{Con}(L)$. We call it the \textbf{determination congruence} of $L$.
For a subset $S \subseteq L$ and congruences $\theta, \psi \in \mathrm{Con}(L)$,  $\theta_{S}$ denotes the restriction of $\theta$ to $S$, and $\theta \circ \psi$ denotes the relational product of $\theta$ and $\psi$.

An algebra $A$ is congruence-permutable if $\theta \circ \psi = \psi \circ \theta$ for all congruences $\theta, \psi$ on $A$. 
As we will see, the congruence $\Phi$ plays a fundamental role in the study of WDLs. The lattice $\mathrm{Con}(L)$ of all congruences of $L$ forms a complete distributive lattice (see \cite{gratz2}).	
\begin{proposition}\label{prop:single}
	Let $\theta$ be a congruence relation on a WDL $L$. The following statements hold:
	\begin{enumerate}
		\item $F = [1]_{\theta}$ (resp. $I = [0]_{\theta}$) is a normal filter (resp. a normal ideal) of $L$. Moreover, $\Psi(I) = F$ and $\Phi(F) = I$.
		
		\item $\alpha = \theta_{\mid S(L)} \in \mathrm{Con}(S(L))$ and $\beta = \theta_{\mid \overline{S}(L)} \in \mathrm{Con}(\overline{S}(L))$.
		
		\item
		\begin{enumerate}
			\item[(i)] If $(0, 1) \in \theta$, then $\theta = L^{2}$.
			\item[(i)] If $S(L) = \overline{S}(L) = \{0, 1\}$, then $\theta = L^{2}$ or $[0]_{\theta} = \{0\}$ and $[1]_{\theta} = \{1\}$.
		\end{enumerate}
		
		\item If $G = [1]_{\theta} \cap \overline{S}(L)$, then $[1]_{\theta} = N[G)$.
		
		\item  $[1]_{\theta} = \{1\}$ if and only if $[0]_{\theta} = \{0\}$.

	\end{enumerate}
\end{proposition}

\begin{proof}
	The first part of (1) is obvious. For the second part, we show that $\Psi(I)=F$, when $I=[0]_{\theta}$ and $F=[1]_{\theta}$. The remaining follows dually. Let $z\in F$, then $(1, z)\in\theta$, by compatibility with $^{\Delta}$ we have $(0, z^{\Delta})\in\theta$, that is $z^{\Delta}\in I$. Since $z\geq z^{\Delta\Delta}$, taking $x=z^{\Delta}$, one has $z\in\Psi(I)$. So $F\subseteq \Psi(I)$. For the reverse inclusion, let $z\in\Psi(I)$, then there exists $u\in I$ such that $z\geq u^{\Delta}$. As $u\in U$, $(u, 0)\in\theta$, so $(u^{\Delta}, 1)\in\theta$, that is $u^{\Delta}\in F$ and $F$ is an upset, so $z\in F$. Hence $F=\Psi(I)$. \\
	Item (2) is straightforward.
	\smallskip
	(3) (i) Suppose that $(0, 1) \in \theta$. Then $1 \in [0]_{\theta}$, and since $[0]_{\theta}$ is an ideal of $L$, it follows that $[0]_{\theta} = L$, hence $\theta = L^{2}$.\\	
	(ii) Assume that $S(L) = \overline{S}(L) = \{0, 1\}$. Then $\theta_{\mid S(L)}$  corresponds to one of the following partitions:	$\{\{0\},\{1\}\}$ or $\{\{0,1\}\}$.
	
	\begin{itemize}
		\item In the first case, $\{\{0\}, \{1\}\}$, we have $[0]_{\theta} = \{0\}$ and $[1]_{\theta} = \{1\}$.
		\item In the second case, $(0, 1) \in \theta$, and by part (3)(i), we obtain $\theta = L^{2}$.
	\end{itemize}
	
	\smallskip
	(4) Let $G = [1]_{\theta} \cap \overline{S}(L)$. We show that $[1]_{\theta} = N[G)$.
	Clearly, $N[G) \subseteq [1]_{\theta}$. For the reverse inclusion, let $z \in [1]_{\theta}$. Then $(z, 1) \in \theta$, hence $(z^{\Delta\Delta}, 1) \in \theta$, and thus $z^{\Delta\Delta} \in G$.
	Since $(z^{\Delta\Delta})^{\Delta\nabla} \leq z$, we deduce by (1) of Proposition~\ref{nor:filter}, that $z \in N[G)$.
	As $z^{\Delta\Delta} \in G$, the desired equality follows.
	
	\smallskip
	(5)
	Assume that $[1]_{\theta} = \{1\}$. Let $a \in [0]_{\theta}$. Then $(a, 0) \in \theta$, and by compatibility with $^{\nabla}$ we have $(a^{\nabla}, 1) \in \theta$.
	By assumption, $a^{\nabla} = 1$, hence $a \leq a^{\nabla\Delta} = 0$, so $a = 0$. Therefore, $[0]_{\theta} = \{0\}$.
	Conversely, assume that $[0]_{\theta} = \{0\}$. Let $a \in [1]_{\theta}$. Then $(a, 1) \in \eta$, and by compatibility of $\theta$ with $^{\Delta}$, we have $(a^{\Delta}, 0) \in \theta$, so $a^{\Delta} = 0$.
	It follows that $a^{\Delta\Delta} = 1 \leq a$, and hence $a = 1$. Therefore, $[1]_{\theta} = \{1\}$.
\end{proof}

\subsection{Congruences generated by  filters}~\vspace{0.2cm}\\
As in the case of Boolean algebras \cite{Burris} and the case of  orthomodular lattices \cite{kalmbach}, we establish a relationship between filters and congruences in distributive WDLs.
Let $F$ be a normal filter of $L$, and consider the binary relation
$$\theta_{F}:=\{(x, y)\in L^{2}\mid \exists u\in F, x\vee u^{\Delta}=y\vee u^{\Delta}\}.$$ 

\begin{theorem}\label{theo: cong induce}
	The binary relation $\theta_{F}$ is the least  congruence relation on $L$ that collapse $F$.
\end{theorem}
\begin{proof}
	(1)  By definition, $\theta_{F}$ is clearly  reflexive and symetric. For transitivity, let $x, y, z\in L$ such that $(x, y), (y, z)\in\theta_{F}$. Then there exist $u, v\in F$ such that $$~ (i) x\vee u^{\Delta}=y\vee u^{\Delta}~\text{ and}~ (ii)~ y\vee v^{\Delta}=z\vee v^{\Delta}.$$ Since $u\wedge v\in F$, using  (i), (ii) and  (4) of Proposition \ref{prop: axiom}, we have $$x\vee (u\wedge v)^{\Delta}=x\vee u^{\Delta}\vee v^{\Delta}=y\vee u^{\Delta}\vee v^{\Delta}=z\vee v^{\Delta}\vee u^{\Delta}=z\vee (u\wedge v)^{\Delta}.$$ 
	Hence $(x, z)\in\theta_{F}$, so $\theta_{F}$ is transitive. Therefore,  $\theta$ is an equivalence relation.\par 
	(2)  We show that $\theta_{F}$ preserves  $\wedge$ and $\vee$.		 
	For $\wedge$. Using distributivity we obtain 	
	$$(x\wedge z)\vee u^{\Delta}=(x\vee u^{\Delta})\wedge (z\vee u^{\Delta})=(y\vee u^{\Delta})\wedge (z\vee u^{\Delta})=(x\wedge z)\vee u^{\Delta}.$$ Hence $\theta_{F}$ is compatible with $\wedge$. 
	
	Let $z\in L$ and assume that $(x, y)\in\theta_{F}$. By commutativity  and  associativity of $\vee$ and $(i)$, $$(x\vee z)\vee u^{\Delta}=(x\vee u^{\Delta})\vee z=(y\vee u^{\Delta})\vee z=(y\vee z)\vee u^{\Delta}.$$ 
	Thus $\theta_{F}$ is compatible with $\vee$.

	\par 
	(3) Next we show that  $\theta_{F}$ is compatible with   the operations $^{\nabla}$ and $^{\Delta}$.\\
	For $^{\nabla}$,  from $x\vee u^{\Delta}=y\vee u^{\Delta}$, by (4') of Proposition \ref{prop: axiom}, we obtain $$x^{\nabla}\wedge u^{\Delta\nabla}=y^{\nabla}\wedge u^{\nabla\Delta}.$$ Since for any $f\in F$, we have   $f\vee f^{\Delta}=1=1\vee f^{\Delta}$, it follows that $(f, 1)\in\theta_{F}$ for any $f\in F$. 
	As $F$ is normal $(u^{\Delta\nabla}, 1)\in \theta_{F}$. Moreover, 
	$(x^{\nabla}\wedge u^{\Delta\nabla},x^{\nabla})\in \theta_{F}$ and $(y^{\nabla}\wedge u^{\Delta\nabla}, y^{\nabla})\in \theta_{F}$. Since   $x^{\nabla}\wedge u^{\Delta\nabla}=y^{\nabla}\wedge u^{\Delta\nabla}$ and $\theta_{F}$ is transitive we deduce that $(x^{\nabla}, y^{\nabla})\in\theta_{F}$. Hence $\theta$ is compatible with $^{\nabla}$.  For $^{\Delta}$, 
	from  $x\vee u^{\Delta}=y\vee u^{\Delta}$, taking the meet  with $u^{\Delta\nabla}$ on both sides give  $$(x\vee u^{\Delta})\wedge u^{\Delta\nabla}=(y\vee u^{\Delta})\wedge u^{\Delta\nabla}$$ that is  $$(x\wedge u^{\Delta\nabla})\vee (u^{\Delta}\wedge u^{\Delta\nabla})=(y\wedge u^{\Delta\nabla})\vee (u^{\Delta\nabla}\wedge u^{\Delta\nabla}).$$
	By (7) of Proposition \ref{prop: axiom} $u^{\Delta}\wedge u^{\Delta\nabla}=0$, hence, 	
	$x\wedge u^{\Delta\nabla}=y\wedge u^{\Delta\nabla}.$  Applying $^{\Delta}$ to both sides yields
	$x^{\Delta}\vee u^{\Delta\nabla\Delta}=y^{\Delta\nabla\Delta}\vee u^{\Delta\nabla\Delta}.$  Since  $u\in F$ and $F$ is normal,  $u^{\Delta\nabla}\in F$ and thus  $(x^{\Delta}, y^{\Delta})\in\theta_{F}$. Thus  $\theta_{F}$ is compatible with $^{\Delta}$. \\
	From (1), (2) and (3) we conclude that   $\theta_{F}$ is a congruence relation of $L$.\par
	\smallskip
	Next we show that $F=[1]_{\theta_{F}}$. Clearly,  $F\subseteq [1]_{\theta_{F}}$. Let $x\in[1]_{\theta_{F}}$. Then $(x, 1)\in\theta_{F}$, so there exists $u\in F$ such that $x\vee u^{\Delta}=1$.  By (7') of Proposition \ref{prop: axiom},  $u^{\Delta\nabla}\leq x$. Since   $u\in F$ and $F$ is normal, $u^{\Delta\nabla}\in F$, hence $x\in F$. Therefore $[1]_{\theta_{F}}=F$. 
	
	\smallskip Finally, we show that $\theta_{F}$ coincide with the congruence generated by $F$. Since $[1]_{\theta_{F}}=F$, it suffices to prove  that if $\psi$ is any congruence on $L$ such that $F\subseteq [1]_{\psi}$, then $\theta_{F}\subseteq \psi$.
	Let $\psi$ be such a  congruence relation, and let $(x, y)\in\theta_{F}$.
	Then there exists $u\in F$ such that $x\vee u^{\Delta}=y\vee u^{\Delta}$. Taking the meet with $u^{\Delta\nabla}$ on both sides give 			 
	$$u^{\Delta\nabla}\wedge (x\vee u^{\Delta})=(x\wedge u^{\Delta\nabla})\vee (u^{\Delta\nabla}\wedge u^{\Delta})=x\wedge u^{\Delta\nabla}=y\wedge u^{\Delta\nabla}\quad (\star_{2}).$$
	Since $(u^{\Delta\nabla}, 1)\in\psi$, and $\Psi$ is  compatible with  $\wedge$, from  ($\star_{2}$) we obtain
	
	$$(x\wedge u^{\Delta\nabla}, x), (y\wedge u^{\Delta\nabla}, y)\in\Psi.$$ 
	By transitivity of $\psi$, it follows that  $(x, y)\in\psi$. Hence $\theta_{F}\subseteq \Psi$, and therefore  $\theta_{F}=\Theta(F)$.
\end{proof}

Assume that $F$ is a normal filter and $(x, y)\in\theta_{F}$. Then,  there exists $u\in F$ such that $x\vee u^{\Delta}=y\vee u^{\Delta}$. Since $F$ is a normal filter, we have   $u^{\Delta\nabla}\in F$. By distributivity and the fact that $x\wedge x^{\nabla}=0$ for all $x\in L$, it follows that  	
$$u^{\Delta\nabla}\wedge (x\vee u^{\Delta})=u^{\Delta\nabla}\wedge x~\text{and}~u^{\Delta\nabla}\wedge (y\vee u^{\Delta})=u^{\Delta\nabla}\wedge y.$$
Hence, $x\wedge u^{\Delta\nabla}=y\wedge u^{\Delta\nabla}$. We therefore obtain the  following remark.

\begin{remark}\label{rek:cons}
	For any $F\in NF(L)$, if $(x, y)\in\theta_{F}$, then there exists $f\in F$ such that $x\wedge f = y\wedge f$.
\end{remark}

Our next result provides a simple rule for computing the join of two congruences in the case where one of them is generated by a filter.

\begin{theorem}\label{theo:congper3}
	Let $F\in NF(L)$ and $\Psi\in Con(L)$. Then~						
	$\theta_{F}\vee \Psi=\theta_{F}\circ\Psi\circ\theta_{F}.$	
	
\end{theorem}
\begin{proof}
	
	Clearly,  we have  $\theta_{F}\circ\Psi\circ\theta_{F}\subseteq \Psi\vee \theta_{F}$.
	For the reverse inclusion,  it suffices to show that $$\Psi\circ\theta_{F}\circ\Psi\subseteq \theta_{F}\circ\Psi\circ \theta_{F}\quad (\ddagger).$$
	Indeed, since  
	$$\Psi\vee \theta_{F}=\bigcup\{\alpha_{i_{1}}\circ\ldots\circ\alpha_{i_{k}}, k\geq 1, \alpha_{i_{l}}\in\{\Psi, \theta_{F}\}, \alpha_{i_{1}}=\Psi \}~\text{(see \cite{Burris})}.$$  The  
	inclusion $(\ddagger)$ implies that
	$$\Psi, \theta_{F}, \Psi\circ\theta_{F}, \theta_{F}\circ\Psi,\Psi\circ\theta_{F}\circ\Psi\circ\alpha_{1}\circ ..\circ \alpha_{k}\subseteq\theta_{F}\circ\Psi\circ\theta_{F},$$ for any~ $\alpha_{j}\in\{\Psi, \theta_{F}\}~\text{and}~, 1\leq j\leq k.$
	Hence the desired equality follows.\par 
	Now, let us prove ($\ddagger$). 
	Let $(a,b)\in\Psi\circ\theta_{F}\circ\Psi$. Then there exist $u, v\in L$ such that $a\Psi u\theta_{F} v\Psi b$.
	Since  $u\theta_{F}v$, by Remark \ref{rek:cons} there exists $t\in F$ such that $u\wedge t=v\wedge t$. It follows that $$a\theta_{F} (a\wedge t), (a\wedge t)\Psi (u\wedge t)=(v\wedge t) \text{and}~ (v\wedge t)\Psi (b\wedge t),  (b\wedge t)\theta_{F} b.$$ Thus,  $$(a, a\wedge t)\in\theta_{F}, (a\wedge t, b\wedge t)\in\Psi~\text{and}~ (b\wedge t, b)\in\theta_{F},$$ which implies that   $(a, b)\in \theta_{F}\circ\Psi\circ\theta_{F}$. Finally, note that   $(b\wedge t, b)\in\theta_{F}$,  because $(b\wedge t)\vee t^{\Delta}=b\vee t^{\Delta}$ by distributivity and $t\in F$. This complete the proof.
\end{proof}	
Recall that an algebra $\mathcal{A}$ is \textbf{regular} if for each $\theta, \beta\in Con(\mathcal{A})$, if there exists $a\in L$ such that $[a]_{\theta}=[a]_{\beta}$, then $\theta=\beta$.
\begin{theorem}\label{theo: max}
	\begin{enumerate}
		\item $\Phi$ is the largest congruence $\eta$ on $L$ such that $[1]_{\eta} = \{1\}$.
		\item A WDL $L$ is regular if and only if $\Phi = \Delta$.
	\end{enumerate}
\end{theorem}

\begin{proof}
	(1) Clearly one has $[1]_{\Phi}=\{1\}$, in fact let $a\in[1]_{\Phi}$, then $a^{\Delta}=a^{\nabla}=0$ and $a^{\Delta\nabla}=1\leq a$,  $a=1$. Thus $[1]_{ \Phi}=\{1\}$. Let $\eta$ be a congruence on $L$ with  $[1]_{\eta}=\{1\}$. Then by (5) of Proposition \ref{prop:single}, one has $[0]_{\eta}=\{0\}$. Let's show that $\eta\subseteq \Phi$. Let $(a, b)\in \eta$. By compatibility with $^{\nabla}$ one has $(a^{\nabla},  b^{\nabla})\in \eta$, therefore;  for any $x\in L$, $(x\vee a^{\nabla},  x\vee b^{\nabla})\in \eta$ $(\star)$.  Taking $x=b^{\nabla\nabla}$ in ($\star$), we obtain  $(b^{\nabla\nabla}\vee b^{\nabla}=1,  b^{\nabla\nabla}\vee a^{\nabla})\in \eta$. By assumption,  $b^{\nabla\nabla}\vee a^{\nabla}=1$ and by (7') of Proposition \ref{prop: axiom} we deduce that $b^{\nabla}\leq a^{\nabla}$. Similarly using $x=a^{\nabla\nabla}$ in $(\star$), one has $a^{\nabla}\leq b^{\nabla}$.  Hence $a^{\nabla}=b^{\nabla}$. Also one has $(a^{\Delta}, b^{\Delta})\in\eta$ implies
	$(x\wedge a^{\Delta}, x\wedge b^{\Delta})\in\eta$ ($\star_{1}$), replacing $x$ by $b^{\Delta\Delta}$ in ($\star_{1}$), and then by $a^{\Delta\Delta}$ and using two times (7) of Proposition \ref{prop: axiom}, one obtains $a^{\Delta}=b^{\Delta}$. Therefore $(a,b)\in\Phi$. Hence $\eta\subseteq\Phi$. \par 
	(2) For the  direct implication, it is clear that if $L$ is regular, then $\Delta=\Phi$, in fact, $[1]_{\Phi}=\{1\}$ and $[0]_{\Phi}=\{0\}$ and $\Delta$ share a  classe with $\Phi$, by regularity one has $\Phi=\Delta$. Now assume that $\Delta=\Phi$ and let's show that $L$ is regular. 
	Now let us consider a congruence $\theta$ on $L$ with $\{b\}$ as a class.  Let $a\in I=[0]_{\theta}$, then $(a, 0)\in\theta$, and $(a^{\nabla\nabla}, 0)\in  \theta$. One deduces that $(b\vee a^{\nabla\nabla},  b)\in \theta$, therefore $a^{\nabla\nabla}\vee b=b$, that is $a^{\nabla\nabla}\leq b$. Similarly $a^{\nabla}\cong 1(\theta)$, implies $a^{\nabla}\wedge b\cong b(\theta)$, that is $a^{\nabla}\wedge b=b$, so $b\leq a^{\nabla}$, it follows that $a\leq a^{\nabla\nabla}\leq a^{\nabla}$, so $a\leq a\wedge a^{\nabla}=0$, hence $a=0$. Thus  $I=[0]_{\theta}=\{0\}$, and $[1]_{\theta}=\{1\}$. By using (1) one has $\theta\subseteq\Phi=\Delta$, that is $\theta=\Delta$.
\end{proof}

One can then define a regular WDL as a WDL that satisfies the condition $\Phi = \Delta_L$.	

\begin{theorem}
	In a distributive WDL, the congruences generated by filters commute.
\end{theorem}

\begin{proof}
	Let $F_{1}$ and $F_{2}$ be two normal filters of $L$. We show that $\theta_{F_{1}} \circ \theta_{F_{2}} \subseteq \theta_{F_{2}} \circ \theta_{F_{1}}$. 
	
	Assume that $x \,\theta_{F_{1}}\, z \,\theta_{F_{2}}\, y$. Then there exist $f_{1} \in F_{1}$ and $f_{2} \in F_{2}$ such that 
	\[
	x \vee f_{1}^{\Delta} = z \vee f_{1}^{\Delta} \quad (\dagger_{1})
	\quad \text{and} \quad
	z \vee f_{2}^{\Delta} = y \vee f_{2}^{\Delta} \quad (\dagger_{2}).
	\]
	
	Let $\omega = (x \vee f_{2}^{\Delta}) \wedge (y \vee f_{1}^{\Delta})$. Note that 
	\[
	(x \vee f_{2}^{\Delta},\, x) \in \theta_{F_{2}} \quad (\dagger_{3})
	\quad \text{and} \quad
	(y \vee f_{1}^{\Delta},\, y) \in \theta_{F_{1}} \quad (\dagger_{4}).
	\]
	
	Since $(y, z) \in \theta_{F_{2}}$, we have $(y \vee f_{1}^{\Delta},\, z \vee f_{1}^{\Delta}) \in \theta_{F_{2}}$. Using $(\dagger_{4})$, it follows that 
	\[
	((y \vee f_{1}^{\Delta}) \wedge (x \vee f_{2}^{\Delta}),\, (z \vee f_{1}^{\Delta}) \wedge x) \in \theta_{F_{2}}.
	\]
	Hence $\omega \,\theta_{F_{2}}\, x$.
	
	Moreover, from $(x, z) \in \theta_{F_{1}}$, we get $(x \vee f_{2}^{\Delta},\, z \vee f_{2}^{\Delta}) \in \theta_{F_{1}}$. Using $(\dagger_{3})$, we obtain 
	\[
	(\omega,\, (z \vee f_{2}^{\Delta}) \wedge y) \in \theta_{F_{1}}.
	\]
	Since $(z \vee f_{2}^{\Delta}) \wedge y = y$, it follows that $(\omega, y) \in \theta_{F_{1}}$ and $(x, \omega) \in \theta_{F_{2}}$, hence $(x, y)\in\theta_{F_{2}}\circ\theta_{F_{1}}$. Therefore, 
	\[
	\theta_{F_{1}} \circ \theta_{F_{2}} \subseteq \theta_{F_{2}} \circ \theta_{F_{1}}.
	\]
	
	A similar argument shows that $\theta_{F_{2}} \circ \theta_{F_{1}} \subseteq \theta_{F_{1}} \circ \theta_{F_{2}}$. Thus, the desired equality holds.
\end{proof}		

\begin{corollary}
	For any $F_{1}, F_{2} \in NF(L)$, we have $\theta_{F_{1}} \vee \theta_{F_{2}} = \theta_{F_{1}} \circ \theta_{F_{2}}$.
\end{corollary}

\begin{proof}
	Setting $\Psi = \theta_{F_{2}}$ and $\theta_{F} = \theta_{F_{1}}$ in Theorem~\ref{theo:congper3} completes the proof.
\end{proof}

Let $h : NF(L) \to Con(L)$ be defined by $h(F) = \theta_{F}$.

\begin{lemma}
	Let $L$ be a distributive WDL, and let $F, G \in NF(L)$. The following equivalence holds:
	\[
	F \subseteq G \ \text{if and only if}\  h(F) \subseteq h(G).
	\]
\end{lemma}

\begin{proof}
	It is straightforward to verify that,  if $F \subseteq G$, then $\theta_{F} \subseteq \theta_{G}$. 
	For the converse, assume that $\theta_{F} \subseteq \theta_{G}$. Then $[1]_{\theta_{F}} \subseteq [1]_{\theta_{G}}$, that is, $F \subseteq G$, since $[1]_{\theta_{F}} = F$ and $[1]_{\theta_{G}}=G$.
\end{proof}

We now present some examples of congruences $\theta$ and their corresponding normal filters (respectively, ideals).

\begin{example}
	Consider the WDL $\mathcal{P}_{6}$.
	\begin{enumerate}
		\item The  subset $F=\{b, 1\}$ of $P_{6}$ is a normal filter of $P_{6}$ and $\theta_{F}$ defined by  its classes is given by $\theta_{F}:=\{\{0, u\}$, $\{b, 1\}\}$. 
		\item Let $\beta$ be the equivalence relation on $\mathcal{P}_{6}$ defined by the classes $\{0, b, v\}$ and $\{u, a, 1\}$.  
		Then $\beta\in Con( \mathcal{P}_{6})$, and $F = \{u, a, 1\} = [1]_{\beta}$ is a normal filter.
	\end{enumerate}
\end{example}

\vspace{0.2cm}

\subsection{Simple and subdirectly irreducible distributive WDLs}~\vspace{0.2cm}\\	
In this section, we characterize simple and subdirectly irreducible distributive WDLs, using congruences generated by filters.

\begin{lemma}
	Let $\theta, \beta \in Con(L)$. The following statements hold:
	
	\begin{enumerate}
		\item $\Phi$ is the largest congruence $\theta$ on $L$ such that its restrictions to $\overline{S}(L)$ and $S(L)$ are both trivial.
		
		\item If $\Phi = \Delta$, then for any $\theta \in Con(L)$, if $[1]_{\theta} = \{1\}$, then $\theta = \Delta_{L}$.
		
		\item (i) $(\theta \vee \beta)_{\mid \overline{S}(L)} = (\theta_{\mid \overline{S}(L)}) \vee (\beta_{\mid \overline{S}(L)})$ and
		(ii) $(\theta \vee \beta)_{\mid S(L)} = (\theta_{\mid S(L)}) \vee (\beta_{\mid S(L)})$.
		
		\item If $\theta$ is a maximal congruence on $L$, then $\Phi \subseteq \theta$.
		
		\item Every lattice congruence contained in $\Phi$ is a congruence. Moreover, if $\theta$ is a congruence relation such that $\Phi \vee \theta = L^{2}$, then $\theta = L^{2}$.
	\end{enumerate}
\end{lemma}

\begin{proof}
	1.  Let $\theta$ be a congruence on $L$ satisfying the stated condition. We show that $\theta \subseteq \Phi$.
	Let $(x, y) \in \theta$. Then $(x^{\Delta\Delta}, y^{\Delta\Delta}) \in \theta_{\overline{S}(L)} = \Delta_{\overline{S}(L)}$ and $(x^{\nabla\nabla}, y^{\nabla\nabla}) \in \theta_{\mid S(L)} = \Delta_{\mid S(L)}$. Hence, $x^{\Delta\Delta} = y^{\Delta\Delta}$ and $x^{\nabla\nabla} = y^{\nabla\nabla}$, which implies $(x, y) \in \Phi$. Therefore, $\theta \subseteq \Phi$.
	
	Item 2. follows from 1. of Theorem~\ref{theo: max}.
	
	3. For (i), it is clear that $(\theta_{\overline{S}(L)}) \vee (\beta_{\overline{S}(L)}) \subseteq (\theta \vee \beta)_{\overline{S}(L)}$.
	For the reverse inclusion, let $a, b \in \overline{S}(L)$ such that $(a, b) \in (\theta \vee \beta)_{\overline{S}(L)}$. Then $(a, b) \in \theta \vee \beta$, and there exist
	$a = c_{0}, c_{1}, \ldots, c_{n-1}, c_{n} = b$ such that $(c_{i}, c_{i+1}) \in \theta \cup \beta$ for $0 \leq i \leq n-1$.
	Therefore, $(c_{i}^{\Delta\Delta}, c_{i+1}^{\Delta\Delta}) \in \theta_{\overline{S}(L)} \cup \beta_{\overline{S}(L)}$ for $0 \leq i \leq n-1$,
	which shows that $(a, b) \in \theta_{\overline{S}(L)} \vee \beta_{\overline{S}(L)}$. Hence, the desired equality holds.
	The equality in (ii) is proved similarly to (i).
	
	4. Assume that $\theta$ is a maximal congruence on $L$. Then either $\theta \vee \Phi = \theta$ or $\theta \vee \Phi = L^{2}$.
	Suppose $\theta \vee \Phi = L^{2}$. Using 3., we obtain
	$(\theta \vee \Phi)_{\mid\overline{S}(L)} = \theta_{\mid \overline{S}(L)} \vee \Delta_{\overline{S}(L)} = \overline{S}(L)^{2}$
	(since $\Phi_{\mid\overline{S}(L)} = \Delta_{\overline{S}(L)}$).
	It follows that $(0, 1) \in \theta$, and thus $\theta = L^{2}$, a contradiction.
	Hence, $\Phi \vee \theta = \theta$, and consequently $\Phi \subseteq \theta$.
	
	5. Assume that $\theta$ is an equivalence relation compatible with $\vee$ and $\wedge$, and that $\theta \subseteq \Phi$.
	Let $(x, y) \in \theta$. Then $x^{\Delta} = y^{\Delta}$ and $x^{\nabla} = y^{\nabla}$.
	Since $\theta$ is reflexive, the result follows.
\end{proof}

\begin{proposition}\label{pro: irre}
	Let $L$ be a distributive WDL such that $\Delta = \Phi$. Then  
	$$\cap\{\theta:   \theta\in Con(L)\setminus\{\Delta_{L}\}\}=\{\theta_{F}: F\in NF(L)\setminus\{\{1\}\}\}.$$  	
\end{proposition}

\begin{proof}
	Recall that $\Delta_{L} = \theta_{\{1\}}$ and $L^{2} = \theta_{L}$.
	Set		
	$$\gamma:=\bigcap \{\theta: \theta \in Con(L)\setminus\{\Delta_{L}\}\},  \text{and}~ \psi:=\bigcap\{ \theta_{F}: F\in  NF(L)\setminus \{\{1\}\}\}.$$ 
	We show that $\gamma = \psi$. Clearly, $\gamma \subseteq \psi$.
	For the reverse inclusion, let $(x, y) \in \psi$. Then, for any normal filter $F \neq \{1\}$ of $L$, we have $(x, y) \in \theta_{F}$ $(\star)$.
	
	Let $\theta$ be a congruence relation on $L$ such that $\theta \neq \Delta_{L}$.
	Since $\Delta_{L} = \Phi$ and $\theta$ is proper, we have $F = [1]_{\theta} \neq \{1\}$,
	so $F$ is a proper normal filter of $L$, and by $(\star)$,
	$(x, y) \in \theta_{F} \subseteq \theta$.
	Hence $(x, y) \in \theta$ for every proper congruence $\theta$ on $L$.
	Thus $\gamma = \psi$.
\end{proof}

The following theorem presents the characterization of simple and subdirectly irreducible regular distributive WDLs.

\begin{theorem}
	If  $\Phi = \Delta_{L}$, then  the following statements hold:
	
	\begin{enumerate}
		\item $L$ is subdirectly irreducible if and only if $NF(L)$ has a unique atom.
		\item $L$ is simple if and only if $\mathrm{Card}(NF(L)) = 2$.
	\end{enumerate}
\end{theorem}

\begin{proof}
	Assume that $\Delta = \Phi$ and that $L$ is subdirectly irreducible. Then there exists a unique atom in $Con(L)$, denoted by $\mu$, with $\mu \neq \Delta_{L}$.
	
	If $[1]_{\mu} = \{1\}$, then $[0]_{\mu} = {0}$ and $\mu \subseteq \Phi = \Delta$ (by 1. of Theorem \ref{theo: max}), which is a contradiction. Hence, set $F = [1]_{\mu}$. Since $\mu$ is an atom, we have $\theta_{F} = \mu$. Therefore, $F$ is an atom in $NF(L)$.
	
	If not, suppose there exists $G \subseteq F$, $G \neq F$, in $NF(L)$. Then $\theta_{G} \subsetneq \theta_{F} \subseteq \mu$, which is a contradiction.
	
	Now, assume that $NF(L)$ has a unique atom, denoted by $F$. By Proposition \ref{pro: irre}, we have	
	$\theta_{F}=\cap\{\theta\in Con(L)\setminus\{\Delta\}\}$,	so $\theta_{F}$ is the unique atom of $Con(L)$, and thus $L$ is subdirectly irreducible.\\
	(2) If $L$ is simple, then
	$\{\Theta_{F}\mid F\in N(F)\}\subseteq Con(L)=\{\Delta_{L},\nabla_{L}\}$,	which implies $\mathrm{Card}(NF(L)) = 2$.
	Conversely, assume that $\mathrm{Card}(NF(L)) = 2$ and that there exists a proper congruence $\phi$ on $L$. Since $\phi$ is proper and $L$ is regular, we must have $L \neq [1]_{\phi} \neq \{1\}$, which is a contradiction. Therefore, $L$ is simple.
\end{proof}

\begin{example}
	\begin{enumerate}
		\item $\mathcal{C}_{2}$ and $\mathcal{C}_{3}$ are simple and subdirectly irreducible WDLs.
		\item $NF(P_{6}) = \{\{1\}, F_{1} = \{b,1\}, F_{2} = \{u, a, 1\}, L\}$ is a four-element Boolean algebra, and there are two atoms in $NF(P_{6})$. Hence, $P_{6}$ is not subdirectly irreducible. It is congruence-permutable.  
		Since $\theta = \theta_{F_{1}}$ and $\Psi = \theta_{F_{2}}$ with $\theta \vee \Psi = \nabla$ and $\theta \cap \Psi = \Delta$, one deduces that $P_{6}$ is isomorphic to $C_{2} \times C_{3}$.
		\item The WDL $L_{6}$ is regular, and $NF(L_{6}) = \{\{1\}, L_{6}\}$. Therefore, $L_{6}$ is simple.
		\item $NF(M_{5}) = \{\Delta_{M_{5}}, \nabla_{M_{5}}\}$, and $M_{5}$ is a regular distributive WDL. Therefore, $M_{5}$ is simple and, consequently, subdirectly irreducible.
	\end{enumerate}
\end{example}

Let $g: Con(L) \to NF(L), \; \theta \mapsto [1]_{\theta}$ and $h: NF(L) \to Con(L), \; F \mapsto \theta_{F}$. Both $h$ and $g$ are well-defined.

\begin{theorem}\label{theo iso filt cong}
	If $L$ is a regular WDL, then $Con(L)$ and $NF(L)$ are isomorphic  lattices.
\end{theorem}
\begin{proof}
	$h$ is an injective order preserving map. 
	Since $L$ is regular, by Proposition \ref{pro: irre} $Con(L) = \{\theta_{F} \mid F \in NF(L)\}$.  
	In fact, for a congruence $\Psi$ in $L$, set $F = [1]_{\Psi}$. Then $\theta_{F}$ is a congruence with $[1]_{\theta_{F}} = F = [1]_{\Psi}$. As $L$ is regular, we have $\theta_{F} = \Psi$. Hence $h$ is surjective. Thus $h$ is an order isomorphism.
\end{proof}

\begin{example}
	In $L_{9}$, $B(L_{9})$ is a four-element Boolean algebra, and $NF(L) = \{\{1\}, \{u, d, 1\}, \{b, a, 1\}, L_{9}\}$ is also a four-element Boolean algebra. A similar result holds for $P_{6}$. In particular, in $L_{9}$, setting $F_{1}=\{u, a, 1\}$ and $F_{2}=\{b, a, 1\}$, one has: 
	$
	F_{1} \cap F_{2} =\{1\}, \quad F_{1} \overline{\vee} F_{2} = L_{9},$
	so $L_{9}$ is isomorphic to $L_{9}/F_{1} \times L_{9}/F_{2}$.
\end{example}

\begin{proposition}
	For any WDL $L$, the quotient algebra $L/\Psi$   is a regular WDL. In particular, $\mathcal{C}_{n}/\Phi$ is isomorphic to $\mathcal{C}_{3}$ for any $n \geq 3$.
\end{proposition}

\begin{proof}
	The first part is clear by the definition of $\Phi$. Note that $\Phi$ is the largest proper congruence on $\mathcal{C}_{n}$ and collapses all elements not in $\{0,1\}$. Hence, 
	$
	\mathcal{C}_{n}/\Phi\cong \mathcal{C}_{3},~n\geq 3.$
\end{proof}

\begin{theorem} 
	A WDL $\mathcal{C}_{n}, n \geq 2$, is subdirectly irreducible if and only if $n \leq 4$.
\end{theorem}

\begin{proof}  
	Note that for any proper congruence $\theta$ in $\mathcal{C}_{n}$, we have $[1]_{\theta} = \{1\}$ and $[0]_{\theta} = \{0\}$.  
	Indeed, assume $a \in \mathcal{C}_{n} \setminus \{0,1\}$ and $(1,a) \in \theta$. Then 
	$
	(1^{\Delta}, a^{\Delta}) = (0,1) \in \theta,
	$
	which implies that $\theta = (\mathcal{C}_{n})^{2}$, a contradiction. Hence, $[1]_{\theta} = \{1\}$. Similarly, one shows that $[0]_{\theta} = \{0\}$.
	
	Now, suppose $n \leq 4$. There are three cases: $n = 2, 3$, or $4$.
	\begin{itemize}
		\item \textbf{Case $n=2$:} This is trivial since $\mathcal{C}_{2}$ is a two-element Boolean algebra, so $\mathcal{C}_{2}$ is simple and therefore subdirectly irreducible.
		\item \textbf{Case $n=3$:} There is no proper congruence on $\mathcal{C}_{3}$, following the observation at the beginning of the proof. Hence, $\mathcal{C}_{3}$ is simple and thus subdirectly irreducible.
		\item \textbf{Case $n=4$:} Let $0 < c_{1} < c_{2} < 1$. The only proper congruence on $\mathcal{C}_{4}$ collapses $c_{1}$ and $c_{2}$. Since $c_{1}^{\Delta} = c_{2}^{\Delta} = 1$ and $c_{1}^{\nabla} = c_{2}^{\nabla} = 0$, this congruence is precisely $\Phi_{4}$. It follows that 
		\[
		Con(\mathcal{C}_{4}) = \{\Delta, \Phi_{4}, \nabla\},
		\] 
		and $\mathcal{C}_{4}$ is subdirectly irreducible because $\Phi_{4}$ is the unique atom of $Con(\mathcal{C}_{4})$.
	\end{itemize}
	
	Now, suppose $n \geq 5$, with $0 < c_{1} < c_{2} \leq c_{3} < \ldots < 1$. Consider the equivalence relations
	$
	\theta_{1} = \Delta_{\mathcal{C}_{n}} \cup \{(c_{1}, c_{2}), (c_{2}, c_{1})\}, \quad 
	\theta_{2} = \Delta_{\mathcal{C}_{n}} \cup \{(c_{2}, c_{3}), (c_{3}, c_{2})\}.
	$
	One can check that $\theta_{1}, \theta_{2} \in Con(\mathcal{C}_{n})$ and 
	$
	\theta_{1} \cap \theta_{2} = \Delta_{\mathcal{C}_{n}}.
	$ 
	Therefore, $\mathcal{C}_{n}$ is not subdirectly irreducible for $n \geq 5$.
\end{proof}	
It is well known (\cite{kwdd}) that a given bounded lattice can be endowed with different structures of weakly dicomplemented lattices.

\begin{definition}
	(\cite{kwdd}, Def. 1.2.2) Let $(^{\Delta_{1}}, ^{\nabla_{1}})$ and $(^{\Delta_{2}}, ^{\nabla_{2}})$ be two weak dicomplementations on $L$. We say that $(^{\Delta_{1}}, ^{\nabla_{1}})$ is \textbf{finer} than $(^{\Delta_{2}}, ^{\nabla_{2}})$ (and write $(^{\Delta_{1}}, ^{\nabla_{1}}) \preceq (^{\Delta_{2}}, ^{\nabla_{2}})$) if 
	\[
	x^{\Delta_{1}} \leq x^{\Delta_{2}} \quad \text{and} \quad x^{\nabla_{1}} \geq x^{\nabla_{2}} \quad \text{for all } x \in L.
	\]
\end{definition}

Sometimes one says that the WDL $(L; \wedge, \vee, ^{\Delta_{1}}, ^{\nabla_{1}}, 0, 1)$ is finer than $(L; \wedge, \vee, ^{\Delta_{2}}, ^{\nabla_{2}}, 0, 1)$. The "finer than" relation is an order relation on the class of all weak dicomplementations on a fixed bounded lattice. It admits a top element, namely the top structure or trivial weak dicomplementation. In the case of a distributive WDL, it also admits a bottom element, namely its double p-algebra structure (\cite{kwdd}).

\begin{proposition}\label{prop: norm}
	\begin{enumerate}
		\item Let $(^{\Delta_{1}}, ^{\nabla_{1}})$ and $(^{\Delta_{2}}, ^{\nabla_{2}})$ be two weak dicomplementations on the bounded lattice $L$. If $(^{\Delta_{1}}, ^{\nabla_{1}})$ is finer than $(^{\Delta_{2}}, ^{\nabla_{2}})$ and $F$ is a normal filter for $(^{\Delta_{2}}, ^{\nabla_{2}})$, then $F$ is normal for $(^{\Delta_{1}}, ^{\nabla_{1}})$. The converse is not always true.
		\item Assume that $(^{\Delta_{1}}, ^{\nabla_{1}})$ is finer than $(^{\Delta_{2}}, ^{\nabla_{2}})$ and that both $(L; \wedge, \vee, ^{\Delta_{1}}, ^{\nabla_{1}}, 0,1)$ and $(L; \wedge, \vee, ^{\Delta_{2}}, ^{\nabla_{2}}, 0,1)$ are regular. If $(L; \wedge, \vee, ^{\Delta_{1}}, ^{\nabla_{1}}, 0,1)$ is simple, then $(L; \wedge, \vee, ^{\Delta_{2}}, ^{\nabla_{2}}, 0,1)$ is also simple.
	\end{enumerate}
\end{proposition}

\begin{proof}
	(1) Assume that $(^{\Delta_{1}}, ^{\nabla_{1}})$ is finer than $(^{\Delta_{2}}, ^{\nabla_{2}})$ and that $F$ is a normal filter for $(^{\Delta_{2}}, ^{\nabla_{2}})$.  
	Let $x \in F$. By assumption, $x^{\Delta_{2}\nabla_{2}} \in F$. Since $x^{\Delta_{1}} \leq x^{\Delta_{2}}$ and $^{\nabla_{2}}$ is order-reversing, we have
	\[
	x^{\Delta_{2}\nabla_{2}} \leq x^{\Delta_{1}\nabla_{2}} \leq x^{\Delta_{1}\nabla_{1}}.
	\]  
	As $F$ is an upset, it follows that $x^{\Delta_{1}\nabla_{1}} \in F$. Therefore, $F$ is normal for $(^{\Delta_{1}}, ^{\nabla_{1}})$.
	
	To show that the converse is not always true, consider the WDL $L_{6}$ and $M_{6}$, with weak dicomplementations $(^{\Delta_{L}}, ^{\nabla_{L}})$ and $(^{\Delta_{M}}, ^{\nabla_{M}})$, respectively. One can verify that 
	\[
	x^{\Delta_{M}} \leq x^{\Delta_{L}} \quad \text{and} \quad x^{\nabla_{M}} \geq x^{\nabla_{L}} \quad \text{for all } x \in L = \{0, u, v, a, b, 1\}.
	\]  
	Hence, $(^{\Delta_{M}}, ^{\nabla_{M}})$ is finer than $(^{\Delta_{L}}, ^{\nabla_{L}})$. Let $F = \{u, a, 1\} \subseteq L$. Then $F$ is normal for $P_{6}$, but not for $L_{6}$, since in $L_{6}$, $u \in F$ and $u^{\Delta_{L}\nabla_{L}} = 0 \notin F$.
	
	(2) Suppose $L$ is simple for $(^{\Delta_{1}}, ^{\nabla_{1}})$. Then there are no proper congruences on $L$ for $(^{\Delta_{1}}, ^{\nabla_{1}})$, and consequently no proper normal filters for $(^{\Delta_{1}}, ^{\nabla_{1}})$. By regularity, we have
	\[
	Con(L; \wedge, \vee, ^{\Delta_{2}}, ^{\nabla_{2}}, 0,1) \subseteq Con(L; \wedge, \vee, ^{\Delta_{1}}, ^{\nabla_{1}}, 0,1) = \{\Delta_{L}, \nabla_{L}\}.
	\]
	Therefore, $(L; \wedge, \vee, ^{\Delta_{2}}, ^{\nabla_{2}}, 0,1)$ is also simple.
\end{proof}	
It is clear that if $\theta$ is a congruence relation on $L$, then the quotient set 
$
L/\theta:= \{[a]_{\theta} \mid a \in L\}
$ 
inherits the structure of a WDL. Moreover, if $\phi$ is a congruence relation on $L/\theta$, then 
$
\phi := \beta/\theta = \{([a]_{\theta}, [b]_{\theta}) \mid (a, b) \in \beta\},$ 
where $\beta$ is a congruence on $L$ containing $\theta$. It follows, by the correspondence theorem (see \cite{Burris}, Theorem 6.20), that each congruence on $L/\theta$ is of this form.

A proper normal filter is called \textbf{purely maximal} if it is a maximal element in $(NF(L), \subseteq)$.
For any normal filter $F$ of $L$, we denote by $L/F$ the  algebra $L/\theta_{F}$.

\begin{proposition} 
	Let $F \in NF(L)$. If $L$ is regular, then the algebra $L/\theta_{F}$ is simple if and only if $F$ is a maximal element in $NF(L)$.
\end{proposition}

\begin{proof}  
	Assume that $L$ is regular and $L/F$ is simple. Then the only congruence on $L$ containing $\theta_{F}$ is $L^{2}$, which implies that there is no proper normal filter on $L$ containing $F$. Hence, $F$ is maximal.
	
	Conversely, assume that $F$ is maximal in $(NF(L), \subseteq)$. By the correspondence theorem, congruences in $L/\theta_{F}$ are of the form $\beta/\theta_{F}$, with $\theta_{F} \subseteq \beta$. Since $L$ is regular, all congruences on $L$ are of the form $\theta_{H}$, with $H \in NF(L)$. As $F$ is maximal in $NF(L)$, it follows that there is no proper congruence on $L/F$. Hence, $L/F$ is simple.
\end{proof}

\begin{theorem}
	Let $\theta$ and $\beta$ be congruence relations on $L$ such that $\theta_{S(L)} = \beta_{S(L)}$ and $\beta_{\overline{S}(L)} = \theta_{\overline{S}(L)}$. Then there exists a natural isomorphism given by:
	\[
	\phi: (L/\theta)/\Phi \longrightarrow (L/\beta)/\Phi, \quad ([x]_{\theta})_{\Phi} \mapsto ([x]_{\beta})_{\Phi}.
	\]
\end{theorem}

\begin{proof}
	First, we show that $\phi$ is well-defined. Let $[y]_{\theta} \in ([x]_{\theta})_{\Phi}$. Then 
	\[
	[x^{\Delta}]_{\theta} = [y^{\Delta}]_{\theta} \quad \text{and} \quad [x^{\nabla}]_{\theta} = [y^{\nabla}]_{\theta}.
	\]  
	Since $\theta_{S(L)} = \beta_{S(L)}$ and $\beta_{\overline{S}(L)} = \theta_{\overline{S}(L)}$, it follows that 
	\[
	[x^{\Delta}]_{\beta} = [y^{\Delta}]_{\beta} \quad \text{and} \quad [x^{\nabla}]_{\beta} = [y^{\nabla}]_{\beta},
	\] 
	which implies $[x]_{\beta} \Phi [y]_{\beta}$. Hence, 
	$
	[x]_{\beta}/\Phi =[y]_{\beta}/\Phi,
	$ 
	proving that $\phi$ is well-defined.
	
	It is straightforward to check that $\phi$ is a homomorphism, using the definitions of the operations on the quotient algebras. Surjectivity follows from 
	$
	\phi^{-1}\Big([x]_{\beta}/\Phi\Big) = [x]_{\theta}/\Phi.
	$ 
	It remains to show that $\phi$ is injective. Suppose 
	$
	[x]_{\theta}/\Phi \neq [y]_{\theta}/\Phi.
	$  
	Then either 
	$
	[x^{\Delta}]_{\theta} \neq [y^{\Delta}]_{\theta} \quad \text{or} \quad [x^{\nabla}]_{\theta} \neq [y^{\nabla}]_{\theta},
	$ 
	which, by the assumptions on $\theta$ and $\beta$, implies 
	$
	([x]_{\beta}, [y]_{\beta}) \notin \Phi.
	$ 
	Therefore, 
	$
	\phi\Big([x]_{\theta}/\Phi\Big) \neq \phi\Big([y]_{\theta}/\Phi\Big),$
	proving that $\phi$ is injective.
	Hence, $\phi$ is an isomorphism.
\end{proof}
Let $X \neq \emptyset$ and $x \in X$. For $S \subseteq L$, define 
\[
L^{X}(x, S) = \{f \in L^{X} \mid f(x) \in S\},
\] 
and for $Y \subseteq L^{2}$, set 
$
L^{X}(x, Y) = \{(f, g) \in (L^{X})^{2} \mid (f(x), g(x)) \in Y\}.
$
Define the maps 

$$
\gamma: F(L) \to F(L^{X}), \quad F \mapsto \gamma(F) = L^{X}(x, F), \text{and}~ 
\mu: Con(L) \to Con(L^{X}), \quad \theta \mapsto L^{X}(x, \theta).$$

\begin{proposition}
	Let $F, G \in F(L)$ and $\theta, \beta \in Con(L)$. The following statements hold:
	\begin{enumerate}
		\item $\gamma$ is an embedding of lattices. In particular, if $F$ is normal, then $\gamma(F)$ is normal. 
		\item $Con(L)$ embeds into $Con(L^{X})$ via $\mu$. In particular, if $\theta = \Theta(F)$, then $\mu(\theta) = \theta_{L^{X}(x, F)}$.
		\item $\{L^{X}(x, F) \mid F \in NF(L)\} \subseteq NF(L^{X})$.
	\end{enumerate}
\end{proposition}

\begin{proof}
	(1) Clearly, $\gamma$ is well-defined and satisfies $F \subseteq G$ if and only if $\gamma(F) \subseteq \gamma(G)$. Hence, $\gamma$ is an injective order-preserving map.  
	If $F$ is normal and $f \in \gamma(F)$, then $f(x) \in F$. Since $F$ is normal, 
	\[
	f^{\Delta\nabla}(x) = f(x)^{\Delta\nabla} \in F,
	\] 
	so $f^{\Delta\nabla} \in \gamma(F)$. Therefore, $\gamma(F)$ is normal.\\
	(3) Follows immediately from (1).
	
	(2) Clearly, $\mu$ is well-defined and preserves the order: $\theta \subseteq \beta$ if and only if $\mu(\theta) \subseteq \mu(\beta)$.  
	
	Indeed, assume $\mu(\theta) \subseteq \mu(\beta)$. Let $(a, b) \in \theta$. Then $(\phi_a, \phi_b) \in L^{X}(x, \theta) \subseteq L^{X}(x, \beta)$, so 
	$
	(\phi_a(x), \phi_b(x)) = (a, b) \in \beta,
	$ 
	which shows that $\theta \subseteq \beta$.  
	
	Conversely, if $\theta \subseteq \beta$, then by definition of $\mu$, we have $\mu(\theta) \subseteq \mu(\beta)$.
\end{proof}

Recall that a class $K$ of algebras is said to have the Congruence Extension Property (CEP) if, for any $A, B\in K$ with $A$ a subalgebra of $B$ and for any congruence $\theta$ of $A$, there exists a congruence $\beta$ on $B$ such that $\beta_{\mid  A}=\theta$.
\begin{theorem}(\textbf{Congruence extension property for distributive WDLs}) Let $M, L$ be two regular and  distributive  WDLs such that $M$ is a subalgebra of $L$.
	If $\theta\in Con(M)$, then there exists $\beta\in Con(L)$ such that $\theta=\beta_{\mid M}$. Moreover, if $L$ is a regular WDL, then $Con(M)$ embeds in $Con(L)$.
\end{theorem}

\begin{proof}
	Let $\theta\in Con(M)$. Since $M$ is  a regular WDL,  there exists $G$ a normal filter of $M$ such that $\theta=\theta_{G}$. Let $F=F_{G}$. By 1. of Lemma \ref{lemmiso}, $F_{G}$ is a normal filter of $L$, and by Theorem \ref{theo: cong induce}, $\theta_{F_{G}}$ is a congruence of $L$. Moreover, one has $[1]_{\theta_{F_{G}}}=F_{G}$ and $F_{G}\cap M=G$. Since $G\subseteq F_{G}$, we deduce that $\theta\subseteq \theta_{F_{G}}$. The map $\phi: Con(\mathcal{M})\to Con(\mathcal{L}), \theta=\theta_{G}\mapsto \phi(\theta)=\theta_{F_{G}}$, is an embedding, by 1. of Theorem \ref{iso: normal} and Theorem \ref{theo iso filt cong}.
\end{proof}

\begin{corollary} Let $M$ and $L$ be two WDLs such that $M$ is a subalgebra of $L$.
	\begin{enumerate}	\item If $L$ is a regular WDL, then $M$ is a regular WDL.
		\item If $L$ is a simple WDL, then any subalgebra of $L$ is simple.
		\item A filter $F$ of a distributive WDL $L$ is normal iff and only if it is a co-kernel of a homomorphism of WDLs.
	\end{enumerate}
\end{corollary}
\begin{proof}
1.	Assume that $L$ is a regular WDL, then by (2) Theorem \ref{theo: max},  $\Phi_{L}=\Delta_{L}$, since $M$ is a subalgebra of $L$, $\Phi_{M}\subseteq \Phi_{L}=\Delta_{L}$, it follows that $\Phi_{M}=\Delta_{M}$ and by Theorem \ref{theo: max} $M$ is a regular WDL.\\
	2. Asuume that $L$ is a simple distributive WDL.Then $Con(L)=\{\Delta_{L}, \nabla_{L})\}$ and by congruence extension property $Con(M)=\{\Delta_{M}, \nabla_{M}\}$.\par 
	3. Assume that $L$ is distributive and  $F$ is a normal filter of $L$. Then $\theta_{F}$ is a congruence of $L$  and $\pi_{F}: L\to L/\theta_{F}, x\mapsto [x]_{\theta_{F}}$ is a homomorphism of WDLs with $co-ker(\pi_{F})=F$. Conversely, if $h: L\to L'$ is a homomtphism of WDL, then $F=\{x\in L\mid f(x)=1\}$ is a normal filter of $L$.
\end{proof}
\begin{theorem}
	(Maximal congruences and semi-simple WDLs)
	\begin{enumerate}
		\item  A congruence $\theta$ is  maximal iff $[1]_{\theta}$ is a maximal normal filter of $L$.
		\item 	A regular  WDL $L$ is semi-simple if $\bigcap\{max(NF(L))\}=\{1\}$.
	\end{enumerate}
\end{theorem}

\section{Conclusion}\label{sec13}

In this work, we first provided a comprehensive list of WDLs and from this list we introduced several subclasses of WDLs. We also characterized normal filters generated by subsets and the join of two normal filters, showing that normal filters form a complete lattice that is not a sublattice of 
$F(L)$. The lattice of filters of the Boolean center of $L$ embeds the lattice of normal filters and this embedding become an isomorphism under some conditions. Moreover, we investigated several properties of normal filters in wdls and proved that the lattice of normal filters is isomorphic to the lattice of normal ideals.
Congruences generated by normal filters were characterized, and additional properties were established in the case of distributive WDLs. Finally, we used congruences generated by normal filters to characterize simple and subdirectly irreducible distributive WDLs, deriving several further results on this class of lattices. Congruence extension property is also established for the classe of distributive WDLs.

In future research, we will focus on constructing Priestley spaces corresponding to distributive WDLs.	
\par \textbf{Credit authorship contribution statement}

\backmatter

\bmhead{Supplementary information}

There is no supplementary file.
\bmhead{Acknowledgements}
The work was partially supported by IMU and Simons Foundation, garnt's N° 640 to visit Bern University of Applied Sciences.
\section*{Declarations}

\textbf{Declaration of competing interest}
The authors have no competing interests.
\section*{Declarations}

\textbf{Declaration of competing interest}
The authors have no competing interests.\\
 \textbf{Funding:} 
The work was partially supported by IMU and Simons Foundation, garnt's N° 640 to visit Bern University of Applied Sciences.\\
\textbf{Conflict of interest/Competing interests:} Not applicable.
\\
\textbf{Consent for publication}\\
 \textbf{Data availability:} Not applicable. \\
 \textbf{Materials availability:} Not applicable.
\noindent
If any of the sections are not relevant to your manuscript, please include the heading and write `Not applicable' for that section. 
	
\end{document}